\documentclass[11pt]{article}
\usepackage{amsmath}
\usepackage{amsthm}
\usepackage{amsfonts}
\usepackage{wrapfig}
\usepackage{url}
\usepackage{upgreek}
\usepackage{verbatim}
\usepackage{tikz}
\usepackage[colorlinks=true,citecolor=blue,urlcolor=blue]{hyperref}
\usepackage[T1]{fontenc}
\usepackage{titling}
\usepackage[hmargin=.720in]{geometry}

\newcommand{\sspcoef}{{\cal{C}}}
\newcommand{\dx}{\Delta x}
\newcommand{\dt}{\Delta t}
\newcommand{\Dt}{\Delta t}
\newcommand{\m}[1]{\mathbf{#1}}
\newcommand{\mD}{\m{D}}
\newcommand{\mA}{\m{A}}
\newcommand{\mAt}{\tilde{\m{A}}}
\newcommand{\mS}{\m{S}}
\newcommand{\mSt}{\tilde{\m{S}}}
\newcommand{\mC}{\m{C}}
\newcommand{\mCt}{\m{\tilde{C}}}

\newcommand{\mI}{\m{I}}
\newcommand{\DtFE}{\dt_{\textup{FE}}}

\newcommand{\ceff}{\sspcoef_{\textup{eff}}}

\newcommand\ve{{\bf e}}
\newcommand\vb{{\bf b}}
\newcommand\vbt{{\tilde{\vb}}}
\newcommand\at{\tilde{a}}
\newcommand\bt{\tilde{b}}

\newcommand\vc{{\bf c}}
\newcommand\vct{{\tilde{\vc}}}
\newcommand{\aij}{\alpha_{i,j}}
\newcommand{\bij}{\beta_{i,j}}

\title{Implicit and Implicit-Explicit Strong Stability Preserving Runge--Kutta Methods with High Linear Order}
\author{Sidafa Conde$^{1}$, 
Sigal Gottlieb$^1$,
Zachary J. Grant$^1$, John N. Shadid$^{2,3}$ \\
{\small $^1$Department of Mathematics, University of Massachusetts, Dartmouth} \\
{\small $^2$ Computational Mathematics Department, Sandia National Laboratories} \\
{\small $^3$Department of Mathematics and Statistics, University of New Mexico}
}
\date{}
\begin{document}

\maketitle
\bibliographystyle{siam}

\vspace{-0.5in}
\abstract{When evolving in time the solution of a hyperbolic partial differential equation, it is often 
desirable to use high order strong stability preserving (SSP)  time discretizations. These time
discretizations preserve the monotonicity properties satisfied by the spatial discretization 
when coupled with the first order forward Euler, under a certain time-step restriction.
While the allowable time-step depends on both the spatial and temporal discretizations, 
 the contribution of the temporal discretization can be isolated by taking the 
ratio of the allowable time-step of the high order method to the 
forward Euler time-step. This ratio is called the {\em strong stability coefficient}.
The search for high order strong stability time-stepping methods with high order and 
large allowable time-step had been an active area of research.  
It is known that implicit  SSP Runge--Kutta methods exist only up to sixth order. 
However, if we restrict ourselves to solving only linear autonomous problems, 
the order conditions simplify and we can find implicit SSP Runge--Kutta methods of any {\em linear order}.  
In the current work we  aim to  find very high linear order implicit SSP Runge--Kutta methods that 
are optimal in terms of allowable time-step.  We then show that if we seek optimal implicit  methods 
with high linear orders (up to $p_{lin} = 9$)  that have nonlinear order $p=3$ or $p=4$, the SSP coefficient is not
 significantly affected, but requiring nonlinear order $p=5$ or $p=6$ does significantly reduce the size of the SSP coefficient.
 We also observe that these implicit methods have SSP coefficients which are up to six times as large as the corresponding 
 explicit methods.
Next, we formulate an optimization problem for implicit-explicit (IMEX) SSP Runge--Kutta methods
and find implicit methods with large linear stability regions that pair with known explicit SSP Runge--Kutta methods
of orders $p_{lin}=3,4,6$ as well as  optimized IMEX SSP Runge--Kutta pairs that have high linear order and nonlinear orders $p=2,3,4$.
These methods are then tested on sample problems to verify order of convergence and to demonstrate the sharpness of the SSP
coefficient and the typical behavior of these methods on test problems.}

\normalsize

\vspace{-0.1in}
\section{Strong Stability Preserving Runge--Kutta methods}
\vspace{-0.05in}
The exact solution of a  hyperbolic conservation law of the form 
\begin{eqnarray}
\label{pde}
	U_t + f(U)_x = 0,
\end{eqnarray}
frequently develops sharp gradients or discontinuities, which may cause significant difficulties
in numerical simulations. For this reason, the development of high order spatial discretizations
that can handle such discontinuities or sharp gradients has been an area of active research for the past 
few decades.  Such methods have special nonlinear non-inner-product stability properties that mimic
some significant physical properties of the exact solution, such as  total variation stability, a maximum principle,
 or positivity. These properties  ensure that when \eqref{pde} is discretized in space, the spatial discretization
 $F(u)$ satisfies
  \begin{equation} \label{eqn:FEcond}
	\| u^n + \dt F(u^{n})  \| \leq \| u^n \|, \quad  0  \leq \dt \leq \Delta t_{FE},
\end{equation}
where $u^n$ is a discrete approximation to $U$ at time $t^n$ and 
 $\| \cdot \|$ is the desired norm, semi-norm, or convex functional. In other words, when the semi-discretized system
\begin{eqnarray}
\label{ode}
u_t = F(u),
\end{eqnarray}
is evolved forward in time using a first order explicit Euler method, the numerical
solution satisfies  the desired strong stability property, as long as the time-step is suitably limited.
In practice we want to use a higher order time integrator instead of the first order forward Euler method
\eqref{eqn:FEcond}, but we still need the strong stability  property 
$	\| u^{n+1} \| \le \|u^n\| $
to be satisfied, perhaps under a modified time-step restriction.

Some higher order Runge--Kutta methods can be decomposed into convex combinations
of forward Euler steps \cite{shu1988b}, so that any  convex functional  property satisfied by \eqref{eqn:FEcond} 
is automatically  {\em preserved}, usually under a different time-step restriction. 
For example, if  the $s$-stage implicit Runge--Kutta method is  written in the form   \cite{shu1988, SSPbook2011},
\begin{eqnarray}
\label{rkSO}
y^{(i)} & = & v_i u^n +  \sum_{j=1}^{s} \left( \aij y^{(j)} + \dt \bij F(y^{(j)}) \right), 
\; \; \; \; i=1, . . ., s \\
 u^{n+1} & = & v_{s+1} u^n +  \sum_{j=1}^{s} \left( \alpha_{s+1,j} y^{(j)} + \dt \beta_{s+1,j} F(y^{(j)}) \right), \nonumber
\end{eqnarray}
it is obvious that each stage can be written as a linear combination of the solution at the 
previous step and forward Euler steps of the stages
\begin{subequations}
\begin{equation}
y^{(i)}  = v_i u^n +  \sum_{j=1}^{s} \aij  \left( y^{(j)} + \dt \frac{\bij}{\aij} F(y^{(j)}) \right),
\end{equation}    
\mbox{and}
\begin{equation}
u^{n+1}  = v_{s+1} u^n +  \sum_{j=1}^{s} \alpha_{s+1,j}  \left( y^{(j)} + \dt \frac{\beta_{s+1,j}}{\alpha_{s+1,j}} F(y^{(j)}) \right).
\end{equation}
\end{subequations}

If  the coefficients $v_i$, $\aij$ and $\bij$ are all non-negative and 
$\bij$ is zero whenever the corresponding  $\aij$ is zero, the stages decompose into a 
convex combination of forward Euler steps of the form \eqref{eqn:FEcond},
with each time step replaced by $ \frac{\bij}{\aij}  \dt$. Each stage is then bounded by
\begin{subequations}
\begin{equation} \label{ccbound}
\| y^{(i)}\|  
\leq v_i \left\|    u^n  \right\|  +  \sum_{j=1}^{s} \aij  \,  \left\| y^{(j)} + \dt \frac{\bij}{\aij} F(y^{(j}) \right\|
\end{equation}
\mbox{with the new solution value}
\begin{equation} \label{ccbound2}
\| u^{n+1} \| 
 \leq   v_{s+1} \left\|    u^n  \right\|  +  \sum_{j=1}^{s}  \alpha_{s+1,j}  \,  \left\| y^{(j)} + \dt \frac{ \beta_{s+1,j} }{ \alpha_{s+1,j} } F(y^{(j}) \right\| .
\end{equation}
\end{subequations}
Recall that $v_i+ \sum_{j=1}^{s} \aij =1$ for consistency, and that 
each $\| y^{(j)} + \dt \frac{\bij}{\aij} F(y^{(j)}) \| \leq \| y^{(j)} \|$ for
$ \frac{\bij}{\aij}  \dt \leq \DtFE$ from  \eqref{eqn:FEcond}, so we have 
\begin{eqnarray}
	\label{rkSSP}
	\| y^{(i)}\| \leq \| u^{n}\|  \; \; \; \mbox{and} \; \; \; \;  \|u^{n+1} \| \leq \| u^n \|  \; \; \; \; \mbox{for any} \; \; \;
	\dt \leq \sspcoef  \DtFE 
\end{eqnarray}
where $\sspcoef = \min_{i,j}  \frac{\aij}{\bij}$.
(If any of the $\bij$'s are equal to zero, we consider the corresponding 
ratios to be infinite.) 
Any method that can be written in this form with $\sspcoef>0$ is called a 
{\em strong stability preserving (SSP)} method.

Strong stability preserving second and third order explicit Runge--Kutta methods \cite{shu1988} 
and fourth order methods \cite{SpiteriRuuth2002,ketcheson2008}
were developed using this convex combination approach. These methods ensure that any 
strong stability property satisfied by the spatial discretization when using the forward Euler condition
 \eqref{eqn:FEcond} is preserved by the higher order strong stability preserving Runge--Kutta method.
Furthermore, the convex combination  decomposition above ensures that the strong stability property is also 
satisfied by  the intermediate stages in a Runge--Kutta method. This may be desirable in many applications,
notably in simulations that require positivity.  Clearly, the condition above is  sufficient  for preservation of strong stability;
in   \cite{ferracina2004, ferracina2005,higueras2004a, higueras2005a} it was  shown that this condition is necessary, as well.

 While the time-step depends on both the spatial and temporal discretizations, 
we isolate the contribution of the temporal discretization to the time-step
restriction by considering the ratio of the allowable time-step of the high order method to the 
forward Euler time-step. This ratio is called the {\em strong stability preserving coefficient}.
Using this approach, we view the time-step restriction \eqref{rkSSP} as  a combination of two
 factors: the forward Euler time-step $\DtFE$ that comes from the spatial discretization, and 
 the SSP coefficient $\sspcoef$ that is a property only of the time-discretization. For this reason,
 the research on SSP methods focuses on optimizing the allowable time-step  $\Dt \le \sspcoef \DtFE$ by maximizing
the   {\em SSP coefficient} $\sspcoef$ of the method. Among methods of similar types, a more 
relevant quantity is the  {\em effective SSP coefficient} $\ceff = \frac{\sspcoef}{s}$, which takes into account 
the computational  cost of the method at each iteration, defined by the number of stages $s$ (typically also the number
of  function evaluations per time-step). However, this measure is not applicable when comparing explicit and implicit methods.
Table \ref{table:RK_ceff} gives SSP coefficients for 
explicit SSP Runge--Kutta methods of orders $p \leq 4$ and implicit methods of orders $p \leq 6$ and different stages, and the
ratio of the SSP methods of corresponding implicit and explicit methods. 
Explicit Runge--Kutta methods (and in fact all explicit general linear methods)
have a bound on the SSP coefficient $\sspcoef \leq s$ \cite{SSPbook2011}, while 
all optimal implicit Runge--Kutta methods have been observed to have $ \sspcoef \ \leq 2s$  \cite{ketcheson2009},
which is only twice that of the explicit method.
However, an implicit methods of a {\em given number} of stages and order typically 
has an SSP coefficient that is greater than  twice that of the corresponding explicit Runge--Kutta method, as shown in 
Table \ref{table:RK_ceff}.

\begin{table}
\begin{center}
{
 \begin{tabular}{|c|ccc|ccccc|ccc|} \hline
 & \multicolumn{3}{c|}{Explicit Methods}  &\multicolumn{5}{c|}{Implicit Methods}  &  \multicolumn{3}{c|}{Ratio of Im/Ex} \\ \hline
$s$ $\backslash$ $p$ & 2 & 3    & 4     & 2 & 3    & 4    & 5    & 6    & 2 & 3    & 4  \\ \hline
1  & -  & -     &  -  					& 2 & -    & -    & -    & -        		&    &     &    \\
2  & 1   & -     &  -  					& 4 & 2.73 & -    & -    & -         		& 4 &  &   \\
3  & 2   & 1 &  -   					& 6 & 4.83 & 2.05 & -    & -     		& 3 & 4.83 &  \\
4  & 3   & 2  &  -  					& 8 & 6.87 & 4.42 &1.14 &        		& 2.67 & 3.44 &  \\
5  & 4   & 2.65 & 1.51  				& 10 & 8.90 & 6.04 & 3.19 &         	& 2.50 & 3.36 & 4.00 \\
6  & 5   & 3.54 & 2.28 				& 12& 10.92 & 7.80 & 4.97 & 0.18  	& 2.40 & 3.08 & 3.42  \\
7  & 6   & 4.27 & 3.29 	   			& 14 & 12.93& 9.19 & 6.21 & 0.26   	& 2.33 & 3.03 & 2.79\\
8  & 7   & 5.12  & 4.15	 			& 16 & 14.94 & 10.67 & 7.56 & 2.25  &2.29 & 2.92 & 2.57 \\
9  & 8   & 6.03 &  4.86	 			& 18 & 16.94& 12.04 & 8.90 & 5.80  & 2.25 & 2.81 & 2.48  \\
10 & 9  & 6.80  & 6.00  				& 20 & 18.95 & 13.64& 10.13 & 8.10 & 2.22 & 2.79 & 2.27  \\
11 & 10 &7.59  & 6.50				& 22 & 20.95 & 15.18 & 11.33 & 8.85 & 2.20 & 2.76 & 2.34  \\
\hline
\end{tabular}}
\caption{SSP coefficients of optimal explicit and implicit Runge--Kutta methods and the ratio of the implicit coefficient to the explicit 
coefficient. 
A dash indicates that SSP methods of this type cannot exist, a blank space indicates none were found.}
\label{table:RK_ceff}
\end{center}
\end{table}

In addition to bounds on the effective SSP coefficients, SSP Runge--Kutta methods also suffer from barriers on their
order:  explicit  Runge--Kutta methods with positive SSP coefficient cannot be more than
fourth-order accurate \cite{kraaijevanger1991,ruuth2001} and implicit Runge--Kutta methods with positive SSP coefficient cannot be more than
sixth-order accurate \cite{ketcheson2009,SSPbook2011}. 
These restrictive order barriers of SSP Runge--Kutta methods is a result of 
the nonlinearity of the ODEs.  However, if we are only interested in the order of accuracy on linear autonomous ODE systems, 
Runge--Kutta methods  need only satisfy a smaller set of  order conditions.
If we only want the method to have  high {\em  linear} order ($p_{lin}$), then the order barrier is 
broken and Runge--Kutta methods with  positive SSP coefficients exist for arbitrarily high linear orders.
These methods are of interest because their SSP coefficients
serve as upper bounds for nonlinear methods. However, they are also useful in their own right for linear problems where the 
strong stability preserving property is required, such as Maxwell's equations and the equations of linear elasticity.

In \cite{LNL} the observation that the order barrier is only applicable to the nonlinear order
led to the study of explicit SSP Runge--Kutta methods with high 
{\em linear order} and optimal  nonlinear order $p=4$, which demonstrated that
where high linear order is desired, going to higher nonlinear order costs little or nothing in terms of the SSP coefficient. 
We refer to methods that may have a higher linear order than nonlinear order  ($p_{lin} \geq p$)
as {\em linear/nonlinear (LNL) methods}.
Table \ref{table:LNL} shows the optimized SSP coefficients of these explicit SSP LNL Runge--Kutta methods.
These methods may be of particular value for problems that require high linear order throughout but also have 
some nonlinear features in some regions which benefit from order  greater than two.

The explicit LNL methods found in \cite{LNL} suggest that implicit methods in this class could be beneficial as well. 
Implicit SSP Runge--Kutta methods with very high linear order have not been widely studied, and these methods could
be very useful for linear problems. In addition, requiring that these methods have a nonlinear order $p>2$ would allow
these methods to be more useful if the problem is generally linear but has regions that feature nonlinearities.
 In this work, we expand upon the explicit methods in \cite{LNL} and seek  optimal implicit SSP methods with very high linear order 
 and nonlinear orders $p\leq 6$.
We first consider implicit SSP Runge--Kutta methods that have high order for linear problems and nonlinear order $p=2$.  
We then consider optimal SSP methods with nonlinear orders $p=3,4,5,6$ and compare the SSP coefficients of these
optimal methods. These methods are presented in Section \ref{Implicit_LNL}. These implicit SSP LNL Runge--Kutta methods
have SSP coefficients that are up to six times the size of the corresponding explicit SSP LNL Runge--Kutta methods.
In Section \ref{Numerical1} we  verify the convergence of these methods as well as the sharpness of the SSP coefficient on sample problems. 
Next, in Section \ref{IMEX_LNL} we consider implicit-explicit  Runge--Kutta methods, for use with problems of the form
\[u_t = F(u) + G(u) \] where we wish to treat $F$ explicitly and $G$ implicitly.
We formulate an optimization problem that depends on the ratio of the forward Euler
condition of the two components $F$ and $G$. Using the optimization routine we developed, we found SSP IMEX Runge--Kutta
methods  where the explicit component has  nonlinear order $p^e=3$ and $p^e=4$ 
and higher linear order $p_{lin}>4$, and the implicit component has nonlinear orders $p^i=2,3,4$  with higher linear order  $p_{lin}$ 
that matches the linear order of the explicit part. 
In Section \ref{IMEX_LNL_methods} we present the optimal methods of this type and 
in Section \ref{Numerical2} we verify the convergence of these methods and  show the behavior of these methods
on sample problems that require the SSP property.

\section{Implicit Runge--Kutta with high linear order} \label{Implicit_LNL}

\subsection{The SSP optimization problem}
Our goal is to find implicit SSP Runge--Kutta methods of a given order and number of stages with the largest possible SSP coefficient $\sspcoef$.
While the most convenient formulation for directly observing the SSP coefficient of a Runge--Kutta method
is the Shu-Osher  form  \eqref{rkSO},  for formulating the optimization problem 
it is more convenient to use  the Butcher form
\begin{eqnarray} \label{butcher}
u^{(i)} & = & u^{n} + \dt \sum_{j=1}^{s} a_{ij} F(u^{(j)}) \; \; \; \; (1 \leq i \leq s) 
\\ \nonumber
u^{n+1} & = & u^n + \dt \sum_{j=1}^s b_j  F(u^{(j)}).
\end{eqnarray}
(where the coefficients $a_{ij}$ are place into the matrix $\mA$ and $b_{j}$ into the row vector $\vb$) \cite{ketcheson2008}. 
The reason this approach is preferable is that the Shu-Osher form of a Runge--Kutta method is not unique,
while the Butcher form is, so that rather than perform a search for an optimal convex combination of the Shu-Osher form,
we seek the  unique optimal Butcher coefficients  \eqref{rkSO}.

Following the approach developed  by Ketcheson \cite{ketcheson2008},and  used in 
\cite{ketcheson2008,  ketcheson2009, Ketcheson2010, SSPbook2011, tsrk, LNL} to find
optimal SSP methods,  we aim to maximize the value of $r$ under the following constraints:
\begin{subequations} \label{optimization}
\begin{align}
 \left( \begin{array}{ll} \mA & 0 \\ \vb & 0  \\ \end{array} \right)
\left(\mI + r 
\left( \begin{array}{ll} \mA & 0 \\ \vb & 0  \\ \end{array} \right)
\right)^{-1} \geq 0   
& \; \; \; \mbox{component wise.}  \\
 \left\| r  \left( \begin{array}{ll} \mA & 0 \\ \vb & 0  \\ \end{array} \right)
\left(\mI + r  \left( \begin{array}{ll} \mA & 0 \\ \vb & 0  \\ \end{array} \right)
\right)^{-1} \right\|_\infty \leq 1 \\
 \tau_k(\mA, b) = 0 \; \; \; \mbox{for} \; \; \;  k=1, . . ., P, \label{tau_conditions}
 \end{align}
 \end{subequations}
 where $\tau_k$ are the order conditions, described below.
This optimization gives the Butcher coefficients $\mA$ and $\vb$ and an optimal value
of the SSP coefficient  $\sspcoef = r$.

The conversion from the  optimal Butcher form to the 
canonical Shu-Osher form is given in  \cite{SSPbook2011}, and given the matrix $\mA$, the vector $\vb$ and the SSP coefficient
$r$, the Shu-Osher coefficients $\alpha$ and $\beta$ can be easily accomplished in MATLAB by 
\begin{verbatim}
    s=size(A,1);
    K=[A;b'];
    G=eye(s)+r*A;
    beta=K/G;
    alpha=r*beta;
\end{verbatim}
where the coefficients $v_i$ are then computed by the consistency condition $v_i + \sum_{j=1}^s \alpha_{i,j} =1$.

Some of the major constraints in the optimization problem come from the 
order conditions  $\tau_k(\mA, \vb)$ above. 
For to demonstrate the correct order of accuracy for nonlinear problems,
it must satisfy the order conditions: 

\smallskip

\noindent $\tau_1(\mA, \vb)$: \; $\vb^T \ve   = 1$ \\
$\tau_2(\mA, \vb)$: \;  $ \vb^T \vc   = \frac{1}{2}$ \\
$\tau_3(\mA, \vb)$: \;    $ \vb^T\mA \vc =  \frac{1}{6}$ \;   and  \; $ \vb^T \vc^2 = \frac{1}{3}$  \\
$\tau_4(\mA, \vb)$: \; $\vb^T\mA^2 \vc = \frac{1}{24}$ \; and \;  $ \vb^T \mA \vc^2  = \frac{1}{12}$ \; and \;
$\vb^T \vc\mA \vc  =  \frac{1}{8}$  \; and \;
$ \vb^T \vc^3  =  \frac{1}{4} $ , \\
where  $\vc= \mA \ve$ and $\ve$ is a vector of ones. Thus, for first order there is only one condition ($P=1$ in Equation
\eqref{tau_conditions} above), 
for second order two ($P=2$), for third order four ($P=4$), and for fourth order eight conditions are needed ($P=8$). 
For fifth order, we require a total of $P=17$ conditions, and  for sixth order $P=37$  conditions must be satisfied.

It is well-known that there are no implicit SSP Runge--Kutta methods greater than sixth order. In this work, we will consider
only implicit methods up to order $p=6$, but we will allow the methods to satisfy higher order $p_{lin}\geq p$ 
when applied to a linear problem.
In this case, the order conditions simplify,  and can be expressed as
\begin{eqnarray}
\tau^{lin}_q(\mA, \vb) = \vb^T \mA^{q-2} \vc = \vb^T \mA^{q-1} \ve = \frac{1}{q!} \; \; \; \; \; \forall q=1, . . . , p_{lin}.
\end{eqnarray}

In the following sections we give the SSP conditions of the optimized methods resulting from the optimization problem given in
\eqref{optimization} with  order conditions 
\[ \tau_k(\mA, \vb)  \; \; \; \mbox{for} \; \; 1 \leq k \leq p \]
and
\[ \tau^{lin}_q(\mA, \vb)  \; \; \;  \mbox{for} \; \; p \leq q \leq p_{lin} \]
which we implemented in a MATLAB code \cite{SSPimplicitLNL_github}. The coefficients of these optimized methods
can be downloaded from \cite{SSPimplicitLNL_github}.


\subsection{Optimal SSP LNL implicit Runge--Kutta methods} \label{implicit_LNL} 
We started our search by considering fully implicit Runge--Kutta methods.
However, as in \cite{LNL} these methods had SSP coefficients that were identical to those of the 
optimized diagonally implicit Runge--Kutta (DIRK) methods, so we proceeded to consider only DIRK methods.  

We first found SSP DIRK methods of up to $s=9$ stages with 
with linear order $2 \leq p_{lin} \leq s+1 $ and nonlinear order $p=2$.  
These methods tend to have nice low-storage forms when written in their Shu-Osher arrays.
Our first observation is that the methods with nonlinear order $p=2$ and linear order $p_{lin}=3$ 
the same  SSP coefficient as the third order ($p=3$) optimal  implicit SSP Runge--Kutta methods listed 
 \cite{ketcheson2009,SSPbook2011}. 
 Clearly, the additional nonlinear order condition imposed did not constrain the methods any further.

\begin{figure}[htb]
    \centering
    \begin{minipage}{0.45\textwidth}
        \centering
        \includegraphics[width=1.025\linewidth]{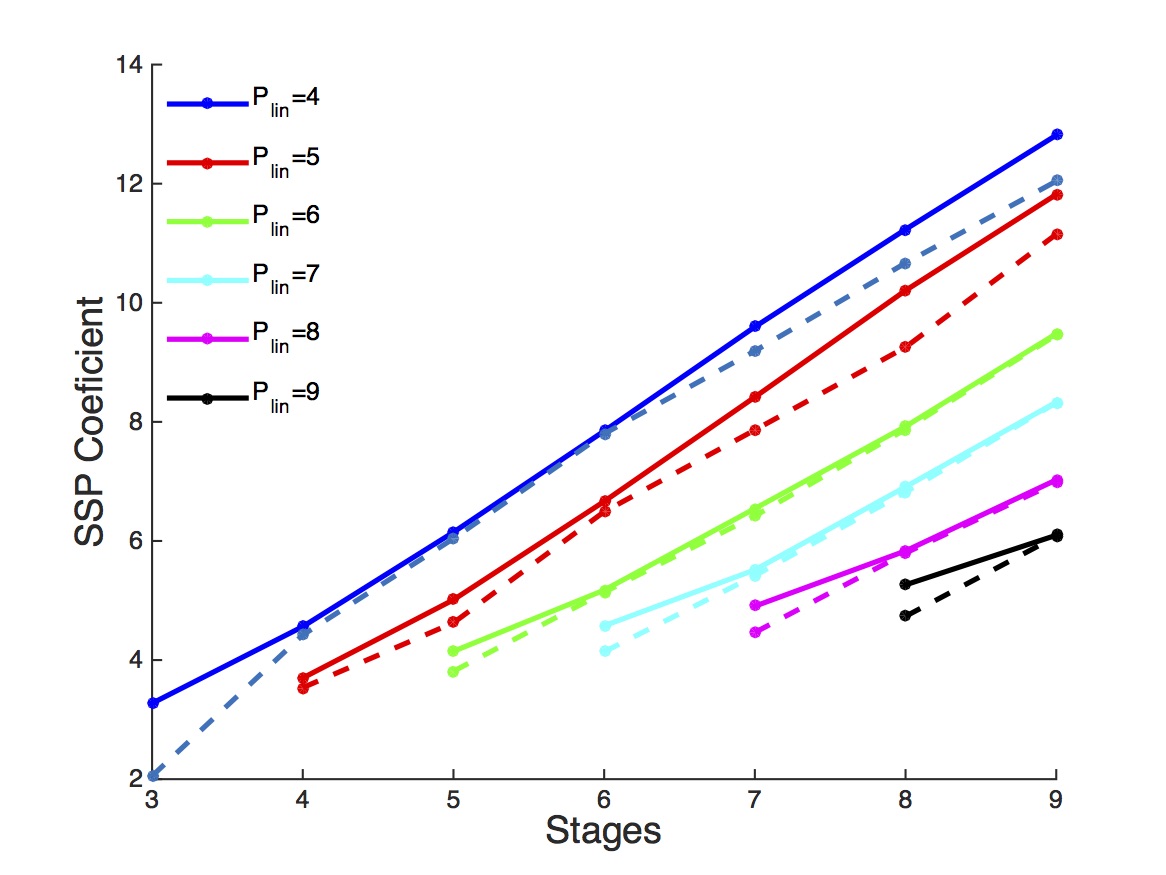}
       \caption{SSP coefficients of methods with nonlinear order $p=2$ (solid lines) and nonlinear order $p=4$ (dashed lines)
       for linear orders $4 \leq p_{lin} \leq 9$, for increasing stages. We note that the SSP coefficients of $p=2$ and $p=3$ 
    are close too identical for all methods with  $s$ stages, and are not 
       significantly different for the two nonlinear orders $p=2$ and $p=4$, and get closer as $s$ increases.}  
\label{fig:SSPiRKcoef1}
    \end{minipage}%
    \hspace{.2in}
    \begin{minipage}{0.45\textwidth}
        \centering
       \includegraphics[width=0.95\linewidth]{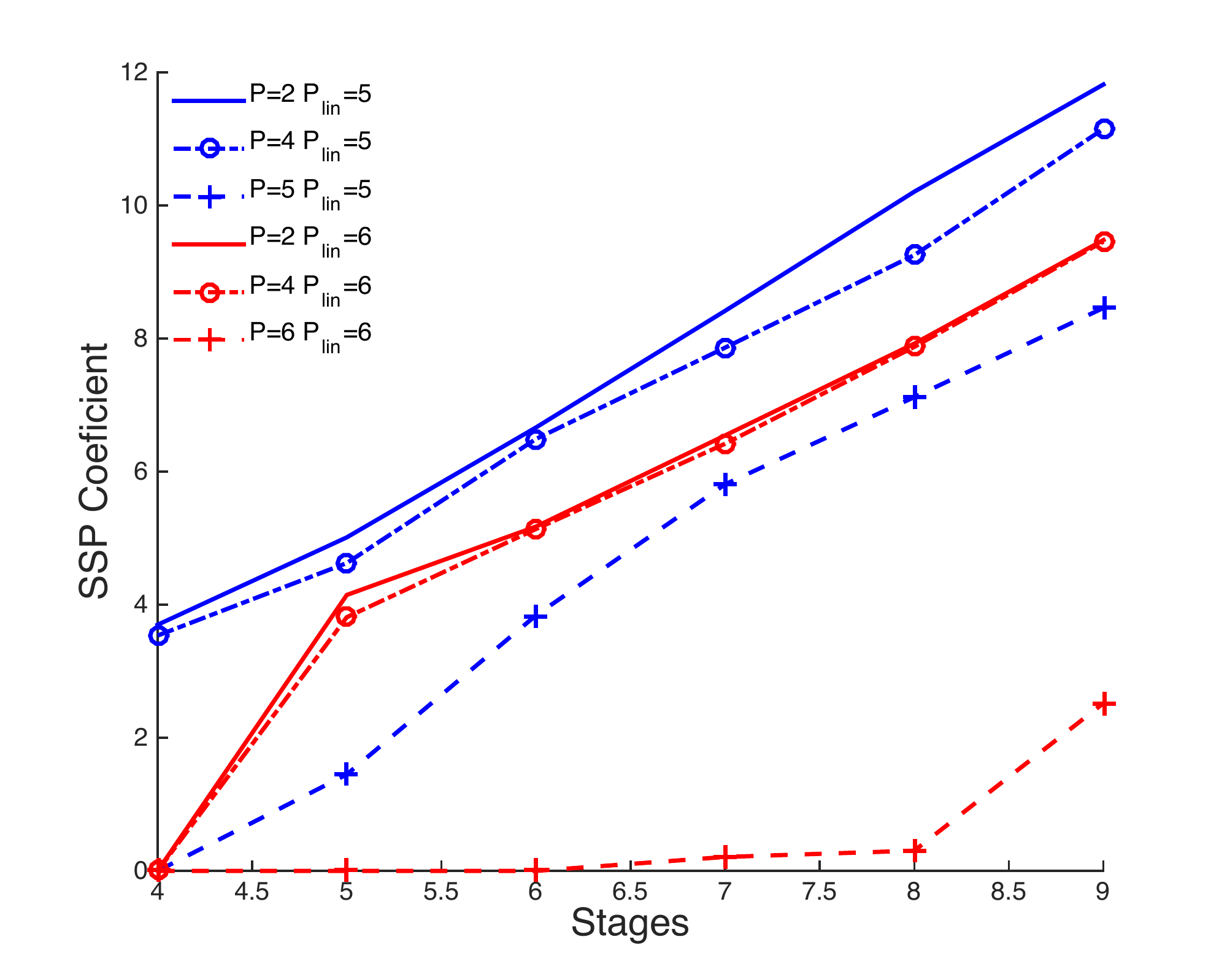}
       \caption{Comparison of SSP coefficients of methods with  linear order $p_{lin}=5$  (blue) and $p_{lin}=6$ (red) for several nonlinear orders.
The $p=2$ lines are solid, $p=4$ lines are dot-dashed with circle marker, and $p=p_{lin}$ lines are  dashed with + markers. 
Clearly, the SSP coefficient of methods with the same number of stages and linear order is significantly smaller for larger nonlinear orders.}
             \label{fig:SSPiRKcoef2} 
           \end{minipage}
\end{figure}

A similar  behavior is observed for SSP DIRK methods with $p_{lin}\geq 4$ and nonlinear orders $p=2$ and $p=3$,
which has SSP coefficients that are the same up to two decimal places,
except for the case of $s=3$ and $p_{lin}=4$ for which there was a difference of $\approx 1.43 \times 10^{-2}$,
and the $s=4$ and $p_{lin}=5$ where there was a difference of $\approx 2 \times 10^{-2}$.  The SSP coefficients are 
 given in Table \ref{tab:ceff} (on left). Finding optimized DIRK SSP Runge--Kutta methods with $p=4$ and 
$\leq 4 p_{lin}\geq 9$, we see that the SSP coefficients in this case are smaller than the methods with nonlinear
orders $p=2,3$, but as the number of stages increases this difference diminishes, as can be seen on 
Table \ref{tab:ceff} (on right). This is also clear in  Figure \ref{fig:SSPiRKcoef1}, which 
plots the SSP coefficients of methods of nonlinear order $p=2$ (and therefore of $p=3$ as well) 
and $p=4$ with different linear orders $3 \leq p_{in} \leq 9$. 
In this figure we observe that if we compare methods with a given $p_{lin}$,
the SSP coefficient of the methods with nonlinear order $p=2$ in  does not have significantly
higher SSP coefficient than the methods with nonlinear order $p=4$. This would seem to suggest that
the linear order is the main constraint on the SSP coefficient, as was the case for 
the explicit methods in \cite{LNL}. However, Figure  \ref{fig:SSPiRKcoef2}, which looks at methods with
higher nonlinear orders $p=5$ and $p=6$ tells a different story.
In this figure we compare three methods  compare three methods with linear order $p_{lin}=5$ and nonlinear orders
$p=2,4,5$ (in blue solid, dot-dashed with circle marker, and dashed with + markers)  and similarly,
three methods with linear order $p_{lin}=6$ and nonlinear orders
$p=2,4,6$ (in red solid, dot-dashed with circle marker, and dashed with + markers)
We say that in this case, as we go to higher nonlinear orders the SSP coefficient is significantly reduced.
For example, a six stage ($s=6$) sixth nonlinear order method has an  SSP coefficient
$\sspcoef=0.18$ whereas the six stage LNL method with nonlinear order $p=4$ and linear order $p_{lin}=6$
has an  SSP coefficient $\sspcoef=5.14$. 

The SSP coefficients of the implicit SSP LNL  Runge--Kutta methods in \ref{tab:ceff} should be compared to those of the
explicit SSP LNL  Runge--Kutta methods \ref{table:LNL}. we see that the implicit methods have SSP coefficients that are
significantly larger than the corresponding explicit methods: up to six times larger. This suggests that the cost of the implicit
solver may be offset by the larger allowable time-step in some cases.
While the methods  in this section are not necessarily suitable for every application,
they are valuable in that they provide an upper bound on the possible SSP coefficient for a given number of
stages and linear and nonlinear orders, and demonstrate the effect of increasing the nonlinear order. 

\begin{table}{\normalsize
\begin{center}
 \begin{tabular}{|l|ccccc||} \hline
&\multicolumn{5}{|c||}{Explicit $p=2,p=3$ methods} \\ \hline
$s \backslash p_{lin}$   &5&6&7&8&9\\\hline
5 &1.00&--&--&--&--\\ \hline
6 &2.00&1.00&--&--&-- \\ \hline
7 &2.65&2&1.00&--&---\\ \hline
8 &3.37&2.65&2.00&1.00&--\\ \hline
9 &4.1&3.37&2.65&2.00&1.00\\ \hline
\end{tabular} \hspace{-.1in}
\begin{tabular}{|ccccc|} \hline
\multicolumn{5}{|c|}{Explicit $p=4$ methods} \\ \hline
5 &6&7&8&9 \\ \hline
0.76 &--&--&--&--\\ \hline
1.81 &0.87&--&--&--\\ \hline
2.57 &1.83 &{\bf 1.00} &--&--\\ \hline
3.36 &2.56&1.93&{\bf 1.00}&--\\ \hline
4.03 &3.35&2.62&1.95&{\bf 1.00}\\ \hline
\end{tabular} 
\caption{SSP coefficients of LNL explicit Runge--Kutta methods from \cite{LNL}.
A dash indicates that SSP methods of this type cannot exist, and the boldface indicates that the SSP coefficients are equal to the 
corresponding ones of lower order.}
\label{table:LNL}
\end{center}}
\end{table}
\begin{table}
\begin{center}
\begin{tabular}{|l|cccccc||cccccc|}
\hline
&\multicolumn{6}{|c||}{Implicit $p=2,p=3$ methods} &\multicolumn{6}{|c|}{Implicit $p=4$ methods}  \\ \hline
$s \backslash p_{lin}$ &4&5&6&7&8&9 &4&5&6&7&8&9\\ \hline
3 &3.24&-- & --&--&--&-- &2.05&-&-&-&-&-  \\ \hline
4 & 4.56&3.64&--&--&--&-- & 4.42&3.54&--&--&--&--  \\ \hline
5 & 6.15&5.01&3.97&--&--&-- &  6.04&4.63 & 3.81&--&--&--  \\ \hline
6 & 7.85&6.66&5.17&4.27&--&-- &  7.80&6.49&5.14&4.14&--&--  \\ \hline
7 & 9.60&8.42&6.54&5.51&4.54&-- & 9.19&7.86&6.42&5.42&4.46&-- \\ \hline
8 & 11.23&10.21&7.92&6.91&5.82 & 4.79 &10.67&9.25&7.86&6.82&5.80&4.74 \\ \hline
9 & 12.81&11.82&9.48&8.33&7.03&6.10
&12.04&11.15&9.46&8.32&6.98&6.08 \\ \hline
\end{tabular}
\caption{\label{tab:ceff} The value of  $\sspcoef$ for LNL iRK with linear order $p_{lin}$, and nonlinear order $p=2$ on left and $p=4$ on right.
Note along the diagonal these the SSP coefficients are close to six times larger than the corresponding explicit LNL methods.}
\end{center}
\end{table}

\subsubsection{Coefficients of selected optimal methods}
The coefficients of all  the methods listed in the section above can be downloaded as {\tt .mat} files from \cite{SSPimplicitLNL_github}. 
In this section we list, for the user's convenience,  the non-zero coefficients of  three methods of particular interest.
We use the canonical Shu--Osher form here, where all the forward-Euler steps are of the same size 
$\frac{\dt}{r}$ and $r=\sspcoef$:
\begin{eqnarray}
\label{cS-Oform}
y^{(i)} & = & v_i u^n +  \sum_{j=1}^{s} \aij  \left( y^{(j)} + \frac{\dt}{r}  F(y^{(j)}) \right), 
\; \; \; \; i=1, . . ., s \\
 u^{n+1} & = & v_{s+1} u^n +  \sum_{j=1}^{s} \alpha_{s+1,j}   \left( y^{(j)} + \frac{dt}{r}  \beta_{s+1,j} F(y^{(j)}) \right), \nonumber
\end{eqnarray}

\noindent {\bf SSP iRK LNL $p=4$, $s=p_{lin}=6$:} One of the most efficient method produced was 
the $s=p_{lin}=6$ LNL method of nonlinear order $p=4$ that has SSP coefficient $\sspcoef =5.138$.
The corresponding implicit Runge--Kutta method in \cite{SSPbook2011} which has   $s=p=p_{lin}=6$ has a much smaller
SSP coefficient of $\sspcoef=0.18$.  The non-zero coefficients are:

\smallskip

{\small
\begin{tabular}{lll}
$\alpha_{1,1}  = 0.227696764527492,$ & $ \alpha_{5,2}=   0.273146312340082, $  & $ \alpha_{7,3}=  0.140604847510042,$  \\
$\alpha_{2,1} = 0.773299008278988, $&  $\alpha_{5,4}= 0.468182990851259,   $& $\alpha_{7,4}= 0.134029552181827,$ \\
$\alpha_{2,2}=   0.226700991721012, $&  $\alpha_{5,5}=  0.226105041192215,  $ & $ \alpha_{7,6}=   0.703213057832428,$ \\
$\alpha_{3,2}=   0.566850708114719, $& $\alpha_{6,5}= 0.772671881656312,    $ & $v_1 = 0.772303235472508,$ \\
$\alpha_{3,3}=  0.245119620891410,  $&  $\alpha_{6,6}=   0.227328118343688,  $ & $v_3 =   0.188029670993872,$ \\
$\alpha_{4,3}=  0.589123375926120,  $&   $\alpha_{7,1}=  0.005835455470528, $ & $v_4 =0.165787716189488,$ \\
$\alpha_{4,4}=   0.245088907884392, $&  $\alpha_{7,2}=  0.016317087005175,  $ & $v_5 = 0.032565655616444,$ \\
\end{tabular}
}
   
\bigskip
   
\noindent {\bf SSP iRK LNL $p=4$, $s=8$, $p_{lin}=9$:} A very high order LNL method with 
$s=8$, $p_{lin}=9$, and $p=4$ has  SSP coefficient $\sspcoef =4.735$. This method has the following non-zero coefficients
given in the canonical Shu-Osher form \eqref{cS-Oform}:

\smallskip

\noindent  {\small    \begin{tabular}{llll}  
$\alpha_{1,1}  =  0.146943975728437,$ &    $\alpha_{5,4}= 0.796548121452431  ,$ &  $\alpha_{7,7} =0.205840405060996     ,$ & 
$\alpha_{9,8}  =0.662855611847356,$ \\
$\alpha_{2,1} = 0.854796464970015   ,$ &$\alpha_{5,5} = 0.136561808924711  ,$ & $\alpha_{8,5}  = 0.510718712707677  ,$ & 
$v_1=  0.853056024271563,$ \\
$\alpha_{2,2} = 0.145203535029985   ,$ & $\alpha_{6,1} =   0.260577803576825  ,$ & $\alpha_{8,7}  =0.353463620808626  ,$ & 
$v_3=    0.251640022410203,$ \\
$\alpha_{3,2} =    0.612204675611763  ,$ & $\alpha_{6,5} = 0.269626835933091  ,$ & $\alpha_{8,8}  =0.135817666483696   ,$ &   
$v_4=    0.122149807578787,$ \\
$\alpha_{3,3} = 0.136155301978034 ,$ & $\alpha_{6,6} =0.206284522717965     ,$ &    $\alpha_{9,1}  =   0.003486997034287  ,$ & 
$v_5=   0.066890069622858,$ \\
$\alpha_{4,3}= 0.742598809241823  ,$ & $\alpha_{7,2} =0.198036604411651 ,$ & $\alpha_{9,2}  =0.067521279383993    ,$ & 
$v_6=  0.263510837772119,$ \\
$\alpha_{4,4}= 0.135251383179389  ,$ &$\alpha_{7,6} =0.596122990527354  ,$ & $\alpha_{9,4}  =0.256478057637965    ,$ &  
$v_9=  0.009658054096400,$ \\ 
    \end{tabular}}

\bigskip

\noindent {\bf SSP iRK LNL $p=2$, $s=10$, $p_{lin}=11$: }  
A very high order LNL method with 
$s=10$, $p_{lin}=11$, and $p=2$ has  SSP coefficient $\sspcoef =5.2306$. This method has a very low storage form, where the  non-zero coefficients
are given in the canonical Shu-Osher form \eqref{cS-Oform}:

\noindent{\small
  \begin{tabular}{llll}   
$\alpha_{1,1}  =    0.193277114534410, $ &    $\alpha_{5,5} =0.117235708890556, $ &   $\alpha_{9,9} = 0.117235191356250, $ &    $v_3 =    0.790541538474273, $ \\          
$\alpha_{2,1} =   0.806723199562524, $ &   $\alpha_{6,5} = 0.718962893859175, $ &   $\alpha_{10,9} = 0.880317745035338, $ &  $v_5 =    0.640785491489029, $ \\ 
$\alpha_{2,2} = 0.193276800437476, $ &     $\alpha_{6,6} =0.117234259419046, $ &       $\alpha_{10,10} =0.117236158521012, $ &      $v_6 =   0.163802846721778, $ \\      
$\alpha_{3,2} = 0.080009844643863, $ &     $\alpha_{7,6} =  0.546025511754727, $ & $\alpha_{11,1} =   0.028409070825259, $ & $v_7 =   0.336736924143727, $ \\
$\alpha_{3,3} =0.129448616881864, $ &     $\alpha_{7,7} =0.117237564101546, $ &    $\alpha_{11,2} = 0.043364313791996, $ &    $v_8 =   0.122161789752829, $ \\ 
$\alpha_{4,3} =   0.870552299752962, $ &    $\alpha_{8,7} =  0.760604303914880, $ & $\alpha_{11,3} = 0.001158601801210, $ &  $v_9 =   0.060130956027420, $ \\
$\alpha_{4,4} = 0.129447700247038, $ &     $\alpha_{8,8} =  0.117233906332291, $ &       $\alpha_{11,10} = 0.921532831100178, $ &     $v_{10} =   0.002446096443650, $ \\      
$\alpha_{5,4} = 0.241978799620415, $ &     $\alpha_{9,8} =   0.822633852616330, $ &  $v_1 = 0.806722885465590, $ & $v_{11} =   0.005535182481357, $ \\
    \end{tabular}
    }
    
    \normalsize
                    
 \subsection{Optimizing for Additional Properties}
\begin{figure}[htb]
    \centering
    \begin{minipage}{0.475\textwidth}
        \centering
        \includegraphics[width=0.975\linewidth]{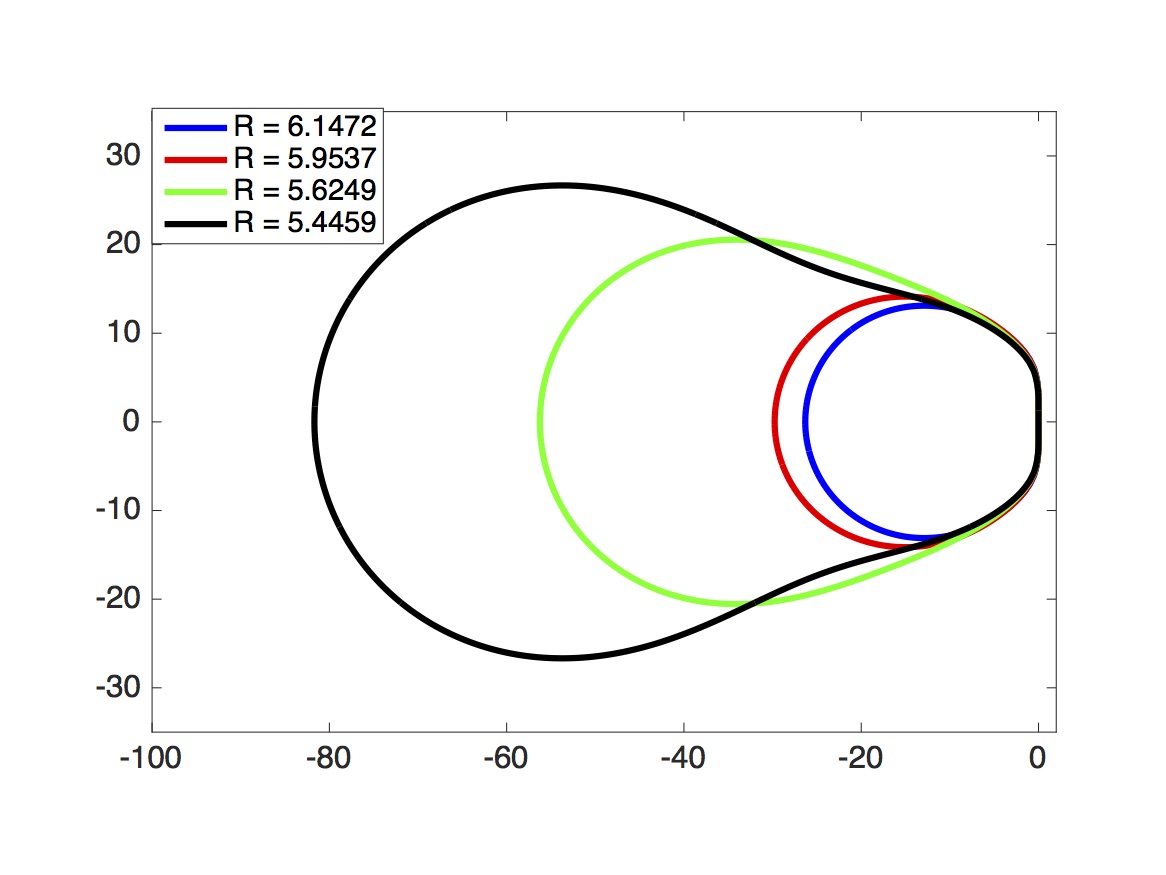}
       \caption{\small Linear stability regions in the complex plane of the optimal SSP method (blue line) and co-optimized methods
       (in red, green, and black) with increasing real axis linear stability.  Shown in the legend are the values of the SSP coefficients for
       these methods.  \vspace{.25in}
        }  
\label{fig:OptNegReal}
    \end{minipage}%
    \hspace{.2in}
    \begin{minipage}{0.475\textwidth}
        \centering
      \includegraphics[width=0.675\linewidth]{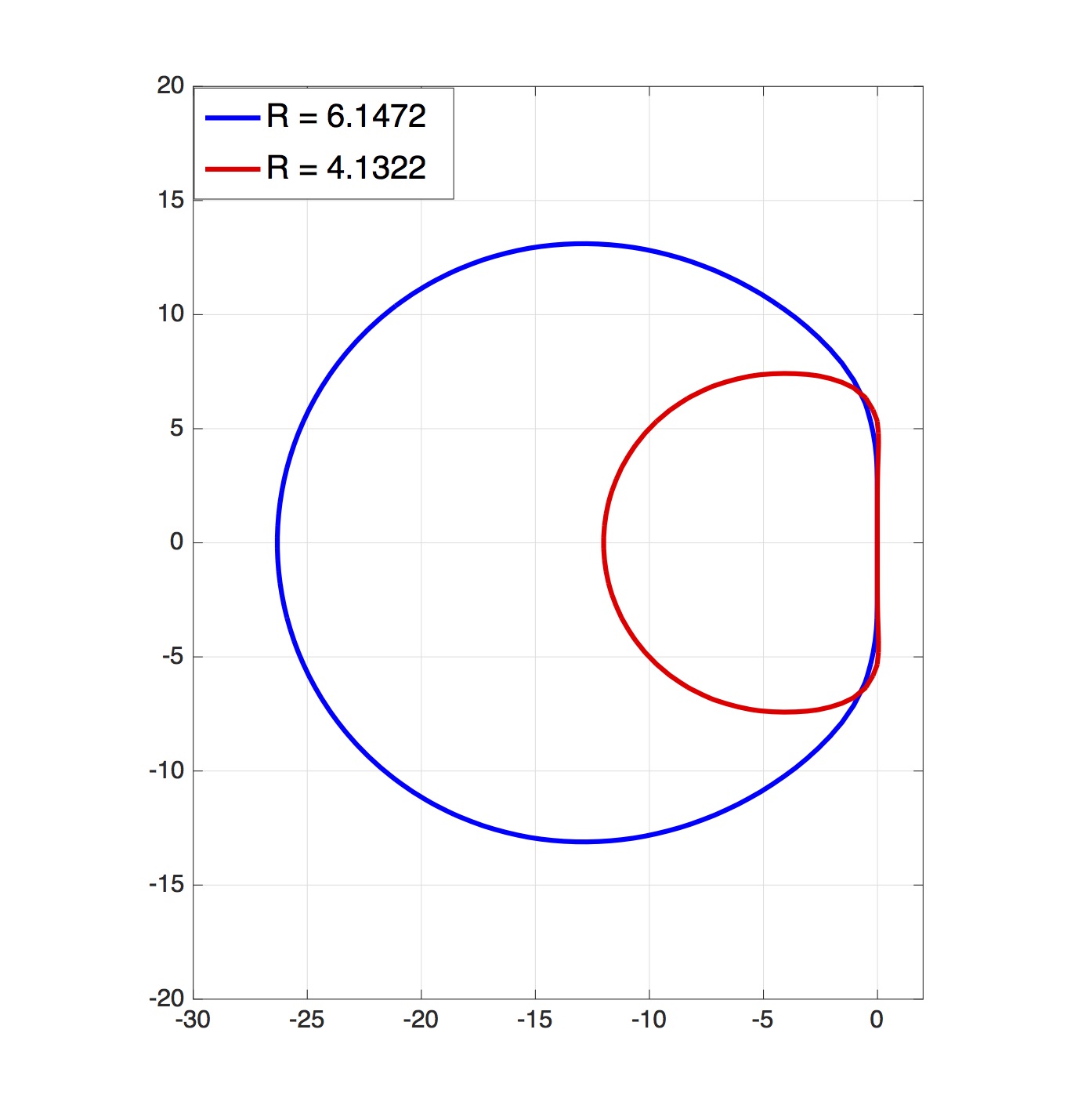} 
        \includegraphics[width=0.2\linewidth, height=.705\linewidth]{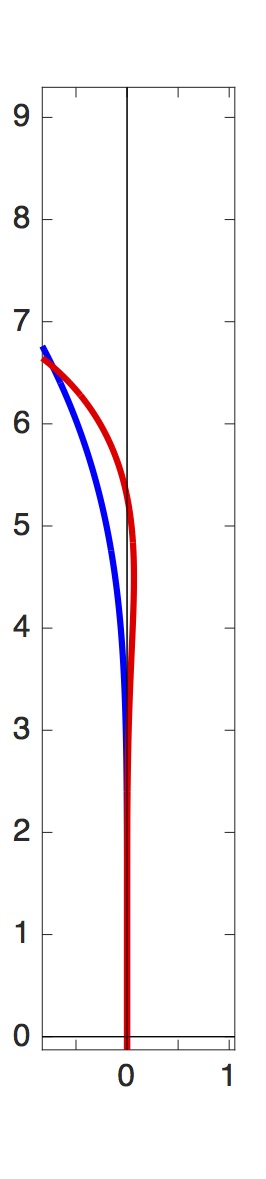}
       \caption{\small  Linear stability regions  in the complex plane of the optimal SSP method (blue line) and co-optimized methods
       in red with larger imaginary axis linear stability, but at the cost of a smaller real axis region and SSP coefficient. On the right is
       the same image zoomed in on the Imaginary axis.
   }
             \label{fig:OptImag} 
           \end{minipage}
\end{figure}

In the section above we optimized the implicit Runge--Kutta methods for the largest possible SSP coefficients. 
However, in many cases the motivation for using an implicit method is not only the SSP coefficient; this is especially true
in cases where the additional time-step allowed for the implicit SSP methods is not enough to offset the cost of 
the implicit solver needed.  If we wish to optimize for other properties, such as linear stability regions,
 alongside the SSP property, we  can add a condition to the inequality constraints in the optimization routine. 
The methods we found above did not have large linear stability regions, and so we wish to explore whether one can find
 methods that have a large SSP coefficient as well as large linear stability region.  In this section we present our optimized methods
 using this approach. The resulting methods are suboptimal in terms of SSP coefficients
but benefit, as desired, from larger regions of absolute stability.  
In Figures \eqref{fig:OptNegReal} and \eqref{fig:OptImag} we show the linear stability regions of methods 
found from this co-optimization approach compared to those found in the subsection above. 

 Figure \eqref{fig:OptNegReal}  shows  several five stage ($s=5$) methods with  $p=3$ and $p_{lin}=4$ resulting from  optimization runs
that required increasing  linear stability regions that included more of the real axis. We show four methods: in blue is the linear stability 
of the optimal SSP method in the subsection above. This method has SSP coefficient $\sspcoef=6.1472$, and crosses the negative real axis
at $x=-26.3$. In red is the linear stability region of a method that has SSP coefficient $\sspcoef= 5.9537$ but allows slightly more of the negative
real axis, crossing it at $x=-29.766$. Next, in green is the linear stability of a method with  a slightly smaller SSP coefficient of $\sspcoef=5.6249$ 
but which allows much more of the negative real axis, crossing it at $x=-56.247$. Finally, in black is the linear stability region of a method with
$\sspcoef=5.4459$ but which allows significantly more of the negative real axis, crossing it at $x=-81.68$. Thus we see that we can co-optimize 
for additional linear stability properties while balancing the need for a large SSP coefficient.

Another frequently desirable feature of time-stepping methods is that they include the imaginary axis or points near the imaginary axis. This is particularly
desirable when solving hyperbolic PDEs. In our case, the optimal SSP five stage ($s=5$) method with  $p=3$ and $p_{lin}=4$ does not include the imaginary axis 
at all ($y=0$). Even when we consider the close neighborhood of the imaginary axis (points which are closer than $10^{-5}$ to the imaginary axis) it only allows these values
up to a value of $y=.9099$. Figure \eqref{fig:OptImag} shows the linear stability of this region in blue.
 If we wish to increase these values, we can obtain a method whose linear stability region is shown in red. This method clearly captures much less
of the real axis, but it includes the imaginary axis up to the value of $y=5.28$. The zoomed image of the region of linear stability in Figure \eqref{fig:OptImag}  
shows this difference  clearly. The SSP coefficient of the second method is significantly less, at $\sspcoef =4.1322$, but this may be a worthwhile trade-off where needed.

The methods in this subsection are meant to represent a small slice of the range of possible properties for which one may co-optimize.
As in \cite{KubatkoKetcheson}, it is possible to optimize SSP methods for particular linear stability regions, 
and it is also possible to consider other desirable properties with respect to which we can to co-optimize.  Furthermore,
other optimization routes are possible, such as starting with a known linear stability polynomial or a particular form of a method that is known to have
desirable properties.

%

\normalsize

\subsection{Numerical Results} \label{Numerical1}
The methods above were found by the optimization code \cite{SSPimplicitLNL_github}, where the order conditions are imposed as equality constraints. It is a
good idea to check that the methods indeed perform as designed on a series of linear and nonlinear test cases. These convergence studies are 
reported in Section \ref{conv1}. Next, we wish to test how the methods perform in terms of the predicted SSP time-step. Clearly, the SSP coefficient
gives a guarantee that the desired property is preserved by the time-stepping method at the corresponding time-step. In Section \ref{tvd1} we study
the difference between the guaranteed time-step and the actual time-step at which the desired nonlinear stability property (in this case the total variation
diminishing property) still holds. A close agreement between these two time-steps demonstrates the relevance of the SSP property for typical cases.

\subsubsection{Verification of the linear and nonlinear orders of convergence} \label{conv1}
\noindent{\bf Example 1.1: Study of the convergence rate for a nonlinear ODE}
To verify the nonlinear order of the implicit Runge--Kutta methods with nonlinear order $p$ and linear order $p_{lin}$
we us a nonlinear system of ODEs
\begin{eqnarray}
&u_1' = u_2  \nonumber \\ 
&u_2' = \frac{1}{\epsilon} (-u_1 + (1-u_1^2) u_2)
\end{eqnarray}
known as the van der Pol problem.
We use $\epsilon = 10$ and initial conditions $u_0 = (0.5;  0)$. 

\begin{figure}[htb]
    \centering
    \begin{minipage}{0.45\textwidth}
        \centering
        \includegraphics[width=0.99\linewidth]{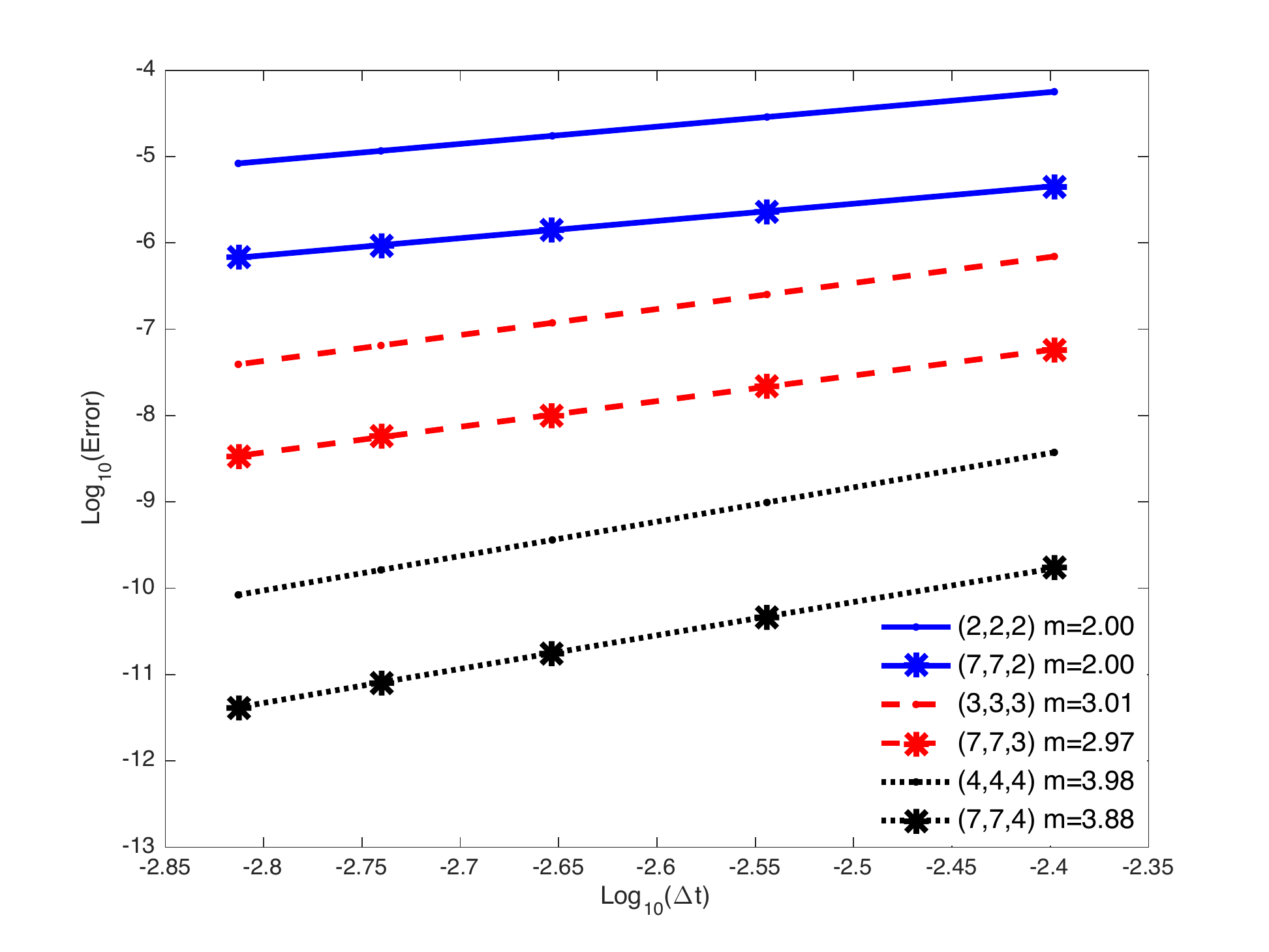}
       \caption{\small Convergence study in  Example 1.1, the nonlinear  van der Pol problem.
   Methods tested are   $(s,p_{lin}, p)= (2,2,2) , (7,7,2) , (3,3,3) ,$ $ (7,7,3) , (4,4,4) , (7,7,4)$, using time-step
$\Dt = \frac{1}{250},  \frac{1}{350},  \frac{1}{450},  \frac{1}{550}, \frac{1}{650}$
to step forward to time $T_{final} = 1.0$. The rate of convergence is, as predicted, the designed
nonlinear order of the method.}  
\label{fig:VDPconv}
    \end{minipage}
    \hspace{.02\textwidth}
    \begin{minipage}{0.45\textwidth}
        \centering
       \includegraphics[width=0.99\linewidth]{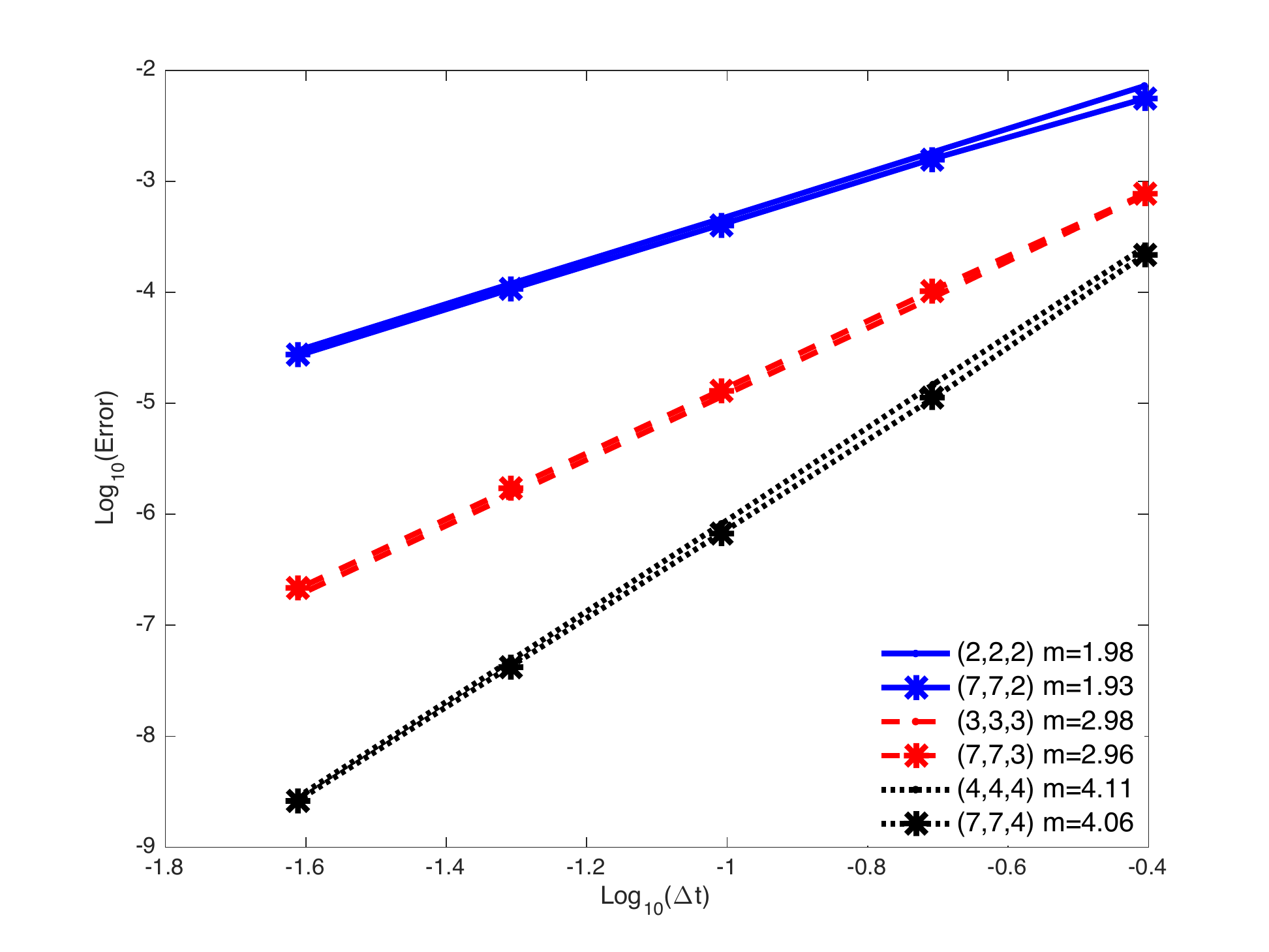}
       \caption{\small Convergence study in Example 1.3, the Buckley-Leverett problem.
       Methods tested are $(s,p_{lin}, p)= (2,2,2) , (7,7,2) ,  (3,3,3) ,$ $ (7,7,3) , (4,4,4) , (7,7,4)$, using time-step
$\Dt = \lambda \Delta x$ with $\lambda = \frac{1}{2}, \frac{1}{4}, \frac{1}{8} , \frac{1}{16},  \frac{1}{32}$, and spatial step
$\dx=\frac{\pi}{4} $ to step forward to time $T_{final} = 2.0$.
The slopes of the lines show the desired order of accuracy. }
             \label{fig:conv_BLspectral}
           \end{minipage}
\end{figure}

This problem was tested with methods with number of stages $s$, linear order $p_{lin}$ and nonlinear order $p$.
Each methods is designated by $(s,p_{lin}, p)$. Of each nonlinear order we test two methods: one with low
stages and $p_{lin} = p$ and one with high number of stages and linear order. 
The methods we test are  $(s,p_{lin}, p)= (2,2,2) , (7,7,2) , (3,3,3) , (7,7,3) , (4,4,4) , (7,7,4)$. 
We  use $\Dt = \frac{1}{250},  \frac{1}{350},  \frac{1}{450},  \frac{1}{550}, \frac{1}{650}$
and step forward to time $T_{final} =  1.0$.
The errors are calculated by using a very accurate solution calculated by  MATLAB's 
ODE45 routine with tolerances set to {\tt AbsTol=}{\tt  RelTol=} $10^{-14}$.
In Figure \ref{fig:VDPconv}  we show that the $log_{10}$ of the errors
in the first component vs. the $log_{10}$ of the time-step $\Dt$. 
The orders (slope of the line) are taken by taking a linear fit using MATLAB's {\tt polyfit}.
We observe the rate of convergence expected for the nonlinear order of the method.

\noindent{\bf Example 1.2: Study of the convergence rate for a linear PDE}
To verify the linear order of  convergence of these methods we solve a linear advection problem $U_t  + U_x=0$
with periodic boundary conditions and initial conditions $U_0(x) =  \sin(x)$ on the spatial domain $x \in [0,2 \pi]$.
We discretize the spatial grid with $N=11$ equidistant points and use the Fourier pseudospectral 
differentiation matrix $\mD$ \cite{HGG2007} to compute $F \approx - U_x$. 
In this case, the solution is a sine wave, so that the  pseudospectral method is exact, and the 
spatial discretization contributes no errors. For this reason, we use  grid refinement in time only,
and use a range of time steps, $\Delta t = \lambda \Delta x$
where  
$\lambda =\frac{1}{10},  \frac{2}{10},  \frac{3}{10},  \frac{4}{10},  \frac{5}{10},  \frac{6}{10},  
\frac{7}{10},  \frac{8}{10},  \frac{9}{10}$
to compute the solution to final time $T_f = 5.0$.
 The errors are measured in the $\ell_2$ norm.

\begin{figure}[htb]
\includegraphics[width=0.475\linewidth]{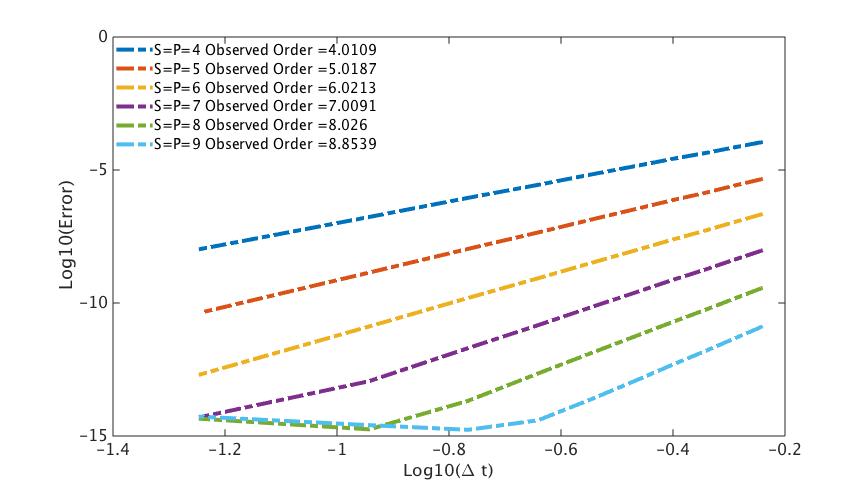}
 \includegraphics[width=.475\textwidth]{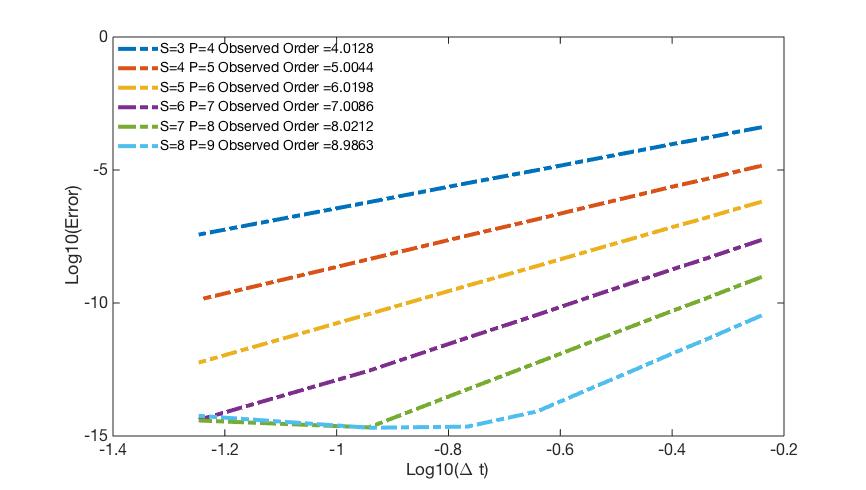}
       \caption{\small  
       Convergence study in  Example 1.2, the linear advection problem with pseudospectral differentiation of the spatial
derivatives. Here we use $N=11$ equidistant points between $(0, 2 \pi)$, and $\Delta t = \lambda \Delta x$. The solution is 
evolved forward to time $T_f=5.0$.
Convergence  plots for linear methods with $s=p_{lin}$ (left) and $s=p_{lin}-1$ (right). The slopes of the lines
       are calculated before round-off error ruins the convergence.}
             \label{fig:conv_spectral}
\end{figure}

Two methods from each order were tested for convergence on this problem.
In Figure \ref{fig:conv_spectral}  we show the convergence plots for the $s=p_{lin} $ methods
the $s=p_{lin} -1$ methods, both with nonlinear order $p=4$.
The slopes were measured in the region before round-off error dominates. 
We observe that the design-order $p_{lin}$ of each method for linear problems is apparent. 

\noindent{\bf Example 1.3: Study of the convergence rate for a nonlinear PDE}
To verify the nonlinear order on PDE, we solve the Buckley-Leverett equation which is 
commonly used to model two-phase flow through porous media:
\begin{align*}
    u_t+f(u)_x & = 0, & \text{ where } f(u) = \frac{u^2}{u^2 +a(1-u)^2},
\end{align*}
on $x\in[0,2\pi]$, with periodic boundary conditions. 
We take $a=\frac{1}{3}$ and initial condition $ u(x,0) = \sin(x) $.
For the spatial discretization, we use a Fourier pseudospectral method as above.
$\Delta t = \lambda \Delta x$ with
$\lambda = \frac{1}{32}, \frac{1}{16}, \frac{1}{8} , \frac{1}{4},  \frac{1}{2}$, and 
$\dx=\frac{\pi}{4} $, evolved to final time $T=2.0$.
The errors are calculated by using a very accurate solution calculated by  MATLAB's 
ODE45 routine with tolerances set to {\tt AbsTol=}{\tt  RelTol=} $10^{-14}$.
In Figure \ref{fig:conv_BLspectral} we show that the $log_{10}$ of the errors
in the first component vs. the $log_{10}$ of the time-step $\Dt$. 
The orders (slope of the line) are taken by taking a linear fit using MATLAB's {\tt polyfit}.
We observe the rate of convergence expected for the nonlinear order of the method.

\subsubsection{Verification of the SSP property} \label{tvd1}

\begin{figure}[t]
\includegraphics[width=0.475\linewidth]{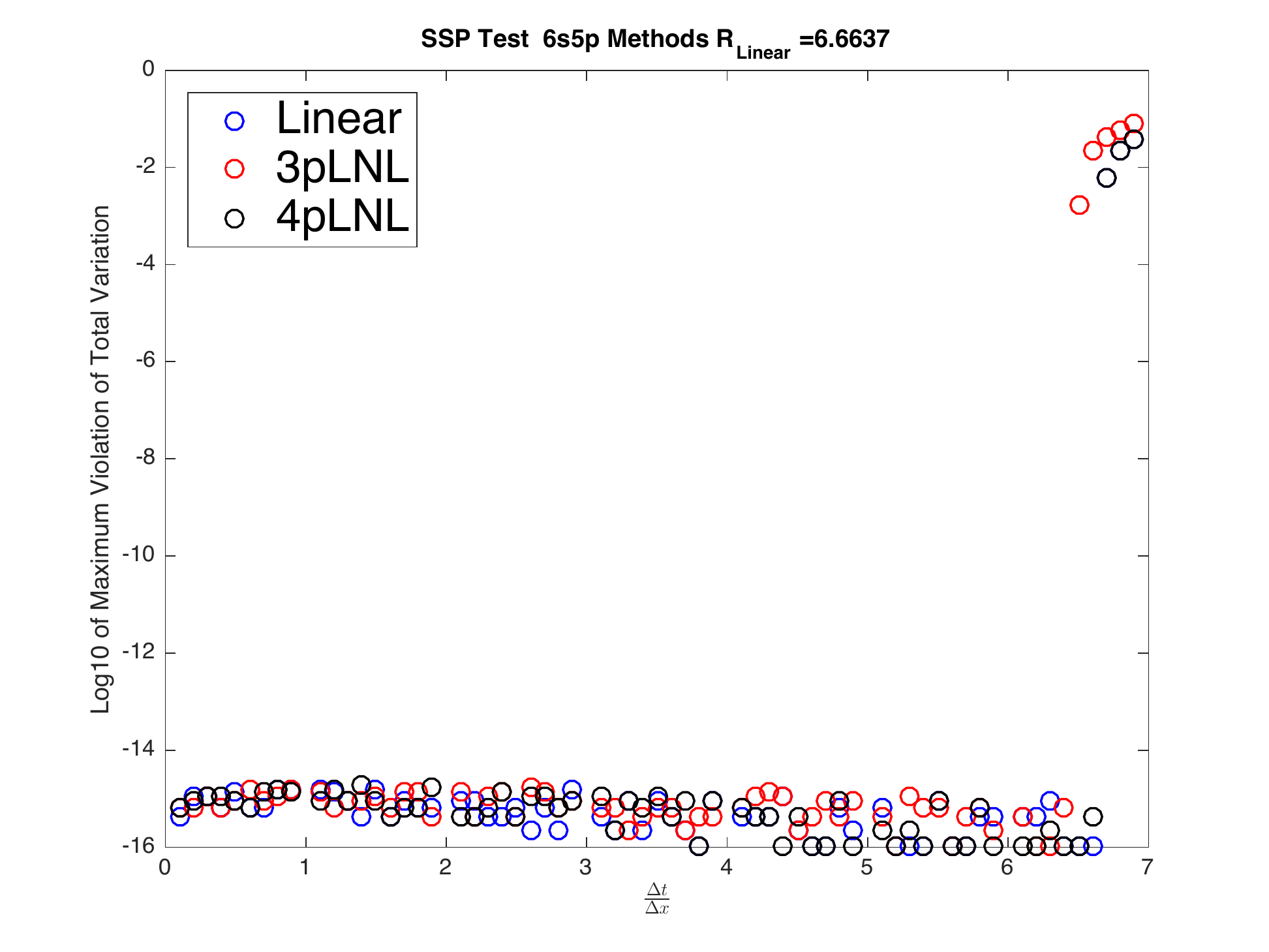}
\includegraphics[width=0.475\linewidth]{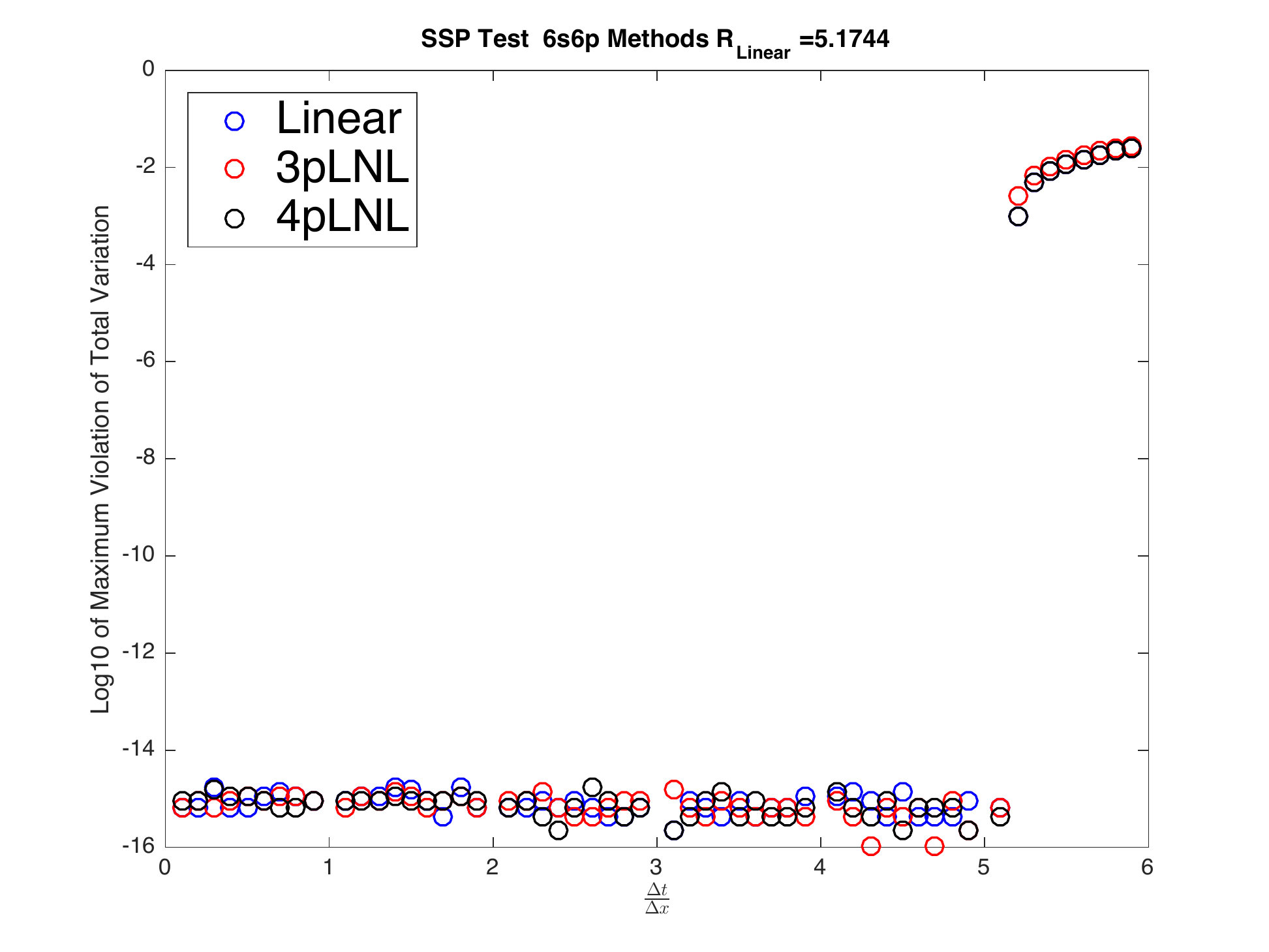}
       \caption{\small  The observed rise in total variation over the first $20$ steps for the $s=6$, $p_{lin}=5$ methods (left) and
       the $s=6$, $p_{lin}=6$ methods. }
       \label{fig:TVtest}
\end{figure}

\noindent{\bf Example 2:} The methods above are selected to optimize the SSP coefficient
$\sspcoef$. This value guarantees that the desired strong stability property that is satisfied by the forward Euler method
will be preserved under the condition $\dt \leq \sspcoef \DtFE$. It is interesting to see
whether for typical problems this value is predictive of the actual time-step at which the desired strong stability 
properties are violated.
As an example, we consider the  linear advection equation
\[ U_t + U_x =0,\]
on the  domain $x \in [-1,1]$ with periodic boundary conditions and the step function initial condition
\[ U_0(x) = \left\{ \begin{array}{rc}
1 &  -0.1 \leq x \leq 0.1 \\
0 & otherwise \\
\end{array} \right. \]
We approximate the spatial derivative with a first order finite difference method, which is total variation diminishing (TVD)
for all values $\dt \leq \dx$.
This simple example is chosen as our experience has shown \cite{SSPbook2011} that this problem often 
demonstrates the sharpness of the SSP time-step.
 For all of our simulations, we use a fixed grid of size $\Delta x = \frac{1}{300}$, and a time-step $\Delta t = \lambda \dx$,
 where we choose values of $\lambda$ starting from the minimum of $\sspcoef/10$ or $0.1$. We then  increase the value
 of $\lambda$ until the TVD property is violated.
  To  measure the effectiveness of these methods, we consider the maximum observed rise in total variation defined by
\begin{equation}
\label{eqn:increase-in-tv}
	\max_{0 \leq n \leq N-1} \left( \| u^{n+1} \|_{TV} - \| u^n \|_{TV} \right). 
\end{equation}
over the first  $N=20$ steps.
We are interested in the time-step in which this rise becomes evident, so we consider a violation if the TV 
rises $10^{-10}$ or more at any time-step. The value of $\lambda$ for which this violation in TV occurs is refined to the level
of  $10^{-12}$.
In Figure \ref{fig:TVtest} we plot the maximal rise in total variation over different values of $\lambda = \frac{\Dt}{\dx}$ 
for a selection of methods. It is interesting to observe that the change in behavior of the methods is quite sharp -- once a 
certain threshold is reached, the SSP coefficient goes from being very small ($<10^{-14}$) 
to being much larger ($> 10^{-4}$) with 
small changes in the time-step. 

We are also interested in difference between the predicted and observed value of $\sspcoef$ (or equivalently, $\dt$) 
at which this violation  happens.   Table \ref{fig:TVD_violationp2} provides these values for the case $p=2$. We observe that
for all the methods with $p=2$  and $p_{lin} < s$, the observed  time-step for which the TVD property is 
preserved matches the predicted time-step up to $ \approx 10^{-10}$.  For these cases, the SSP coefficient
is an excellent predictor of the actual value for which the TVD property breaks down.


\bigskip

\begin{table}[h]
{\small 
\begin{tabular}{|l|l|l|l|l|l|l|l|l|}
\hline
 $s \backslash p_{lin} $&2&3&4&5&6&7&8&9 \\\hline
2& 2.5$\times 10^{-12}$&1.5$\times 10^{-12}$ &&&&&&\\\hline
3&4.0$ \times 10^{-12}$&2.5$ \times 10^{-12}$&3.8$\times 10^{-2}$&&&&&\\\hline
4&8.1$ \times 10^{-12}$&4.5$ \times 10^{-12}$&2.2$\times 10^{-11}$&6.4$\times10^{-2}$&&&&\\\hline
5&1.6$ \times 10^{-11}$&9.1$ \times 10^{-12}$&1.2$\times 10^{-11}$&7.6$ \times 10^{-12}$
&1.7 $\times10^{-1}$ &&&\\\hline
6&3.2$ \times 10^{-11}$&1.8$ \times 10^{-11}$&1.1$\times 10^{-11}$&6.6$ \times 10^{-12}$&2.6$ \times 10^{-11}$
&3.0$\times 10^{-1}$&&\\\hline
7&6.4$ \times 10^{-11}$&3.7$ \times 10^{-11}$&1.4$\times 10^{-11}$&7.6$ \times 10^{-12}$
&2.6$ \times 10^{-10}$&7.1$ \times 10^{-12}$&3.6$\times10^{-1}$&\\\hline
8&1.2$ \times 10^{-10}$&7.4$ \times 10^{-11}$&3.7$\times 10^{-10}$&1.0$ \times 10^{-11}$
&5.5$ \times 10^{-11}$&5.5$ \times 10^{-11}$&1.4$ \times 10^{-2}$&4.7$\times10^{-1}$\\\hline
9&2.5$ \times 10^{-10}$&1.5$ \times 10^{-10}$&1.4$\times 10^{-10}$&1.7$ \times 10^{-10}$&3.4$ \times 10^{-11}$&
2.1$ \times 10^{-11}$&2.2$ \times 10^{-10}$& 3.4$ \times 10^{-2}$ \\ \hline
\end{tabular}
}
\caption{The difference between the theoretical and observed SSP coefficient at which the maximal rise in total
variation is about the threshold of $10^{-10}$ for the $p=2$ methods.}
\label{fig:TVD_violationp2}
\end{table}

For the cases where $p_{lin} = s \leq 7$, and $p=2$,  the observed time step is also within  $\approx 10^{-10}$ of the predicted value.
This is true also for the methods with $p=3$.
However, when $p_{lin} = s =8,9$, and $p=2,3,4$,  the observed time step is much larger than the predicted time-step
 ($\approx 10^{-1}$ or $10^{-2}$). Also, for $p=4$, the $s=p_{lin}=5$, we also have a larger difference between the observed and 
 predicted time-step.
Finally, in the cases where $s = p_{lin}-1$, the observed time step is usually 
much larger than the predicted time-step ($\approx 10^{-1}$ or $10^{-2}$),
and the SSP coefficient is not sharp at all.

\section{Optimal SSP IMEX Runge--Kutta methods with $p_{lin} \geq p$}  \label{IMEX_LNL}
For cases in which the problem has a linear component which restricts the time-step significantly, and a nonlinear component which does not,
 it may be advantageous to treat the two components differently by using an implicit-explicit method. 
 Consider an  initial value problem of the form
\begin{equation}
U'(t) = f(U(t)) + g(U(t)), \qquad U(t_0) = U_0 \qquad t \geq t_0
\label{eq:SC_additiveODE}
\end{equation}
which is semi-discretized to give an ODE of the form 
\begin{equation}
u'(t) = F(u) + G(u).
\label{eq:SC_semi_discrteAdditiveODE}
\end{equation}
To step the two components forward in time, one explicitly and one implicitly, 
we use an Implicit-Explicit (IMEX) Runge--Kutta method, which is a particular kind of  additive Runge--Kutta method.
An $s$-stage Additive Runge--Kutta (ARK) method is defined by two $s \times s$ 
real matrices $\mA, \mAt$ and two real vectors $\vb, \vbt$ such that:
        \vspace{0.1em}
        \begin{equation}
        \begin{aligned}
        u^{(i)} &= u^n + \dt \sum_{j=1}^s a_{ij} F(u^{(j)} + \dt \sum_{j=1}^s \at_{ij} G(u^{(j)} )   \qquad i= 1, \dots,s  \\
                u^{n+1} &= u^n + \dt \sum_{i=1}^s b_i F(u^{(i)} + h \sum_{i = 1}^s \bt_i G(u^{(i)}), \qquad i = 1, \dots, s .
                \label{eq:SC_ark_butcher_form}
                \end{aligned}
                \end{equation}
We select $\mA$ to be a strictly lower triangular matrix, so that this additive Runge--Kutta 
method contains an explicit method (represented by ($\mA, \vb$)) used for the non-stiff part $F$.
As above, we limit ourselves to diagonally implicit methods  and require 
$\mAt$ to be a lower triangular matrix, where ($\mAt, \vbt$) represents an implicit Runge--Kutta  (DIRK)
method used for the stiff part $G$, as in   \cite{Higueras2009,Ascher1997IMEX_ARK, CarpenterKennedyIMEX2003}. 

IMEX methods  were first introduced by Crouziex in 1980 \cite{CrouziexIMEX} for evolving parabolic equations. In 1997, Ascher, Ruuth, and Wetton \cite{Ascher1997IMEX_ms} introduced IMEX multi-step methods for time dependent PDEs, notably convection-diffusion equations, 
and in the same year Ascher, Ruuth and Spiteri  \cite{Ascher1997IMEX_ARK} presented the IMEX Runge--Kutta schemes for such problems.
Although the authors were not focused on designing SSP pairs, some of methods are in fact known SSP methods. 
For example, the method in \cite[Section 2.4]{Ascher1997IMEX_ARK} is the midpoint explicit and implicit SSP methods. 
In this work, all the implicit methods are SDIRK methods, and most have nice properties such as L-stability, with diagonal values
$\gamma$ that are those mentioned in the books of Hairer, Norsett, and Wanner \cite{Hairer1, Hairer2}.

Implicit methods are often particularly desirable when applied to a linear component. In this case, the order conditions simplify.
In \cite{calvo2001}, Calvo, Frutos, and Novo developed IMEX pairs for the case where the implicit  component $G$ is linear. 
This was the first work where the linear and nonlinear orders were separated. The methods they produced had nonlinear implicit order
 $p_{im} = 2$, nonlinear explicit order and linear order $ p_{ex} = p_{lin}=3,4$.

Kennedy and Carpenter \cite{CarpenterKennedyIMEX2003} 
derived IMEX Runge--Kutta methods based on singly diagonally implicit 
Runge--Kutta (SDIRK) methods. This work introduced sophisticated IMEX methods with good accuracy and stability properties,
as well as high quality embedded methods for error control and other features that make these methods usable in complicated 
applications.

The first IMEX methods that had SSP properties were considered by Pareschi and Russo in  \cite{PareschiRusso2005}.
In this work, explicit SSP  Runge--Kutta methods  L-stable implicit components. These schemes were designed to be asymptotically preserving. 
The authors listed a $p_{ex} = p_{im} = p = 3$ order conditions, and also presented a table depicting the number of coupling conditions under
 each underlying assumptions. In this work, the authors designed their methods while enforcing the condition that the abscissas of the explicit
 and implicit methods were equal $\vc = \tilde{\vc}$.  They observed that under the assumption $\vc=\vct$ and $\vb=\vbt$ the 
coupling conditions were redundant for the $p = 3$.   It is worth noting the SDIRK implicit pairs from \cite{PareschiRusso2005} have the same 
$\gamma$ value listed in the books of Hairer, Norsett, and Wanner  \cite{Hairer1, Hairer2}. Further work by Boscarino, Pareschi, and Russo
\cite{BoscarinoPareschiRusso2013} observed that the order reduction phenomenon disappears when $\vb = \vbt$.  In this work 
they presented  a globally stiffly accurate, third order method with non-negative explicit coefficients (named BPR(3,5,3).

Higueras \cite{Higueras2006} was first to consider SSP IMEX methods where both explicit and implicit pairs of the methods were SSP.
In this work, she presented conditions for an IMEX Runge--Kutta methods to be SSP, and listed some previously known methods,
 such as methods from \cite{Ascher1997IMEX_ARK}, that satisfy the SSP conditions. The methods that appear in this work are
 of  order $p \leq 3$, and there is no distinction made between the explicit, implicit and linear order of the methods.
 In  \cite{Higueras2009} Higueras presented order barriers and other characteristics of strong stability 
 preserving additive Runge--Kutta methods. This work contains several very important properties of SSP IMEX pairs
 and is necessary to the understanding of the structure of these methods.
  In later work, \cite{Higueras_TVD_IMEX2012} several known methods  were analyzed in terms
 of their linear stability region and SSP performance on astrophysics-type problems. The implicit second order methods considered
 are L-stable (and thus A-stable as well), but the  third order method was not A-stable. 
 In their work in \cite{higuerasOptSSPIMEXRK2014}, the authors mainly focused on second order, three-stage explicit methods that 
implicit SSP pairs. In this work, the main focus was on methods with {\em suboptimal} SSP coefficients that offered significant 
additional benefits (including large linear stability regions).  

Our approach to deriving optimal SSP methods is slightly different from  \cite{Higueras2006,higuerasOptSSPIMEXRK2014}.
First, we use a slightly different formulation of the SSP condition, which facilitates the construction of 
an optimization problem. Next, we introduce a coefficient which relates the strong stability condition of $F$ to
that of $G$, and make an observation about  time step of optimal methods using this coefficient,
which also simplifies the optimization routine. But the main distinction between our work and prior work is the 
fact that our investigations focus on higher order methods, and especially methods with higher linear order.

To analyze an IMEX method in the SSP context, we assume as in \cite{Higueras2006,higuerasOptSSPIMEXRK2014}
 that when each  component ($F$ or $G$) 
of the ODE is stepped forward individually with a forward Euler method, the new solution 
will satisfy a strong stability property in some desired convex functional $\| \cdot \|$,
but under very different time-step restrictions:
\begin{eqnarray}
	\| u^n + \dt F(u^{n})  \| \leq \| u^n \|, \quad  0  \leq \dt \leq \Delta t_{FE},
 \label{eq:F_fe_condition}
\end{eqnarray}
and
\begin{eqnarray}
	\| u^n + \dt G(u^{n})  \| \leq \| u^n \|, \quad  0  \leq \dt \leq K \Delta t_{FE}.
 \label{eq:G_fe_condition}
\end{eqnarray}
If $K$ is small, the $G$ component will make the overall allowable time-step of the method very small,
so it may be worthwhile to treat this component implicitly. We know \cite{ketcheson2009, SSPbook2011} that the allowable time-step
for an $s$-stage SSP implicit Runge--Kutta method of order $p>1$ can only be expected to be, at most, $2s$ times the allowable 
explicit forward Euler time-step $\DtFE$, which does not typically offset the additional cost of solving the implicit system. 
However, in the previous discussions we show that the allowable time-step of an implicit method is in fact more than twice, 
and in the new methods presented in Section \ref{Implicit_LNL} up to six times, 
as large as that of an explicit method with the same order and 
number  of stages and order.
Furthermore, in the case where $G$ is a linear operator, the cost of the implicit solver is much smaller and under such circumstances 
it may be worthwhile to use an implicit SSP Runge--Kutta method for this part of the problem. However, if 
$F$ is nonlinear, or if for any other reason the cost of implicitly solving the system for $F$ is large,
we wish to use an explicit time-stepping method for this part.
The IMEX approach may also be desirable if the value of $K$ is not small (perhaps even infinitely large)
but other factors, such as linear stability requirements, greatly limit the 
 step size if $G$ is treated explicitly. 
 
The main contribution of our work is to present SSP IMEX methods with  $p_{lin} \geq 4$. These methods have not appeared
 in the literature, and while it is not clear yet how useful they will be in actual applications, we believe that 
 understanding the SSP bounds on such methods and its dependence on the linear and nonlinear orders is of value.

\subsection{Formulating the Optimization Problem}

To formulate the optimization method for this problem, we stack the matrix $\mA$ and the vector $\vb$, padded with zeroes, 
into a square matrix $\mS$  (as we did in Equation \eqref{optimization}). Similarly, we convert $\mAt$ and $\vbt$ into a matrix $\mSt$.
We now rewrite the  method   \eqref{eq:SC_ark_butcher_form} as
\begin{eqnarray} \label{MSMDvector}
Y = \ve u^n + \Delta t \mS  F(Y) +  \Delta t \mSt G(Y).
\end{eqnarray}
Now, we add the terms 
\begin{eqnarray*}
\left( I +r \mS  + \tilde{r} \mSt \right)  Y &=& 
\ve u^n  + r \mS \left( Y + \frac{\Delta t}{r}  F(Y) \right) + \tilde{r} \mSt  \left( Y +  \frac{\Delta t}{ \tilde{r}} G(Y) \right),   \\
 Y & = & R (\ve u^n) + P \left( Y + \frac{\Delta t}{r}  F(Y) \right) + Q  \left( Y +  \frac{\Delta t}{ \tilde{r}} G(Y) \right),
\end{eqnarray*}
where \[
	R  =  \left( I +r \mS  + \tilde{r} \mSt \right)^{-1}, \; \; \; \; \; \;
	P  =  r R \mS,  \; \; \; \; \; \;
	Q =   \tilde{r}  R \mSt. 
\]
From this formulation we see that if $R \ve$, $P$, and $Q$ are all positive component-wise, then the resulting method
is simply a convex combinations of forward Euler steps $Y + \frac{\Delta t}{r}  F(Y)$ and $Y +  \frac{\Delta t}{ \tilde{r}} G(Y)$,
and therefore will preserve the strong stability property in the desired convex functional $\| \cdot \|$,
under the time-step restriction $\dt \leq \min \left(r \DtFE, \tilde{r} K \DtFE \right) $. We observe that the optimal methods will have
$r = \tilde{r} K$, so the optimization problem becomes: \\
{\em Maximize $r$ such that}
\begin{subequations}
\begin{align}
 \left( I +r \mS  + \frac{r}{K}  \mSt \right)^{-1} \ve   \geq 0     \label{Opt_ARK1} \\
r  \left( I +r \mS  + \frac{r}{K}  \mSt \right)^{-1} \mS \geq 0 \label{Opt_ARK2} \\
\frac{r}{K}  \left( I +r \mS  + \frac{r}{K}  \mSt \right)^{-1} \mSt \geq 0  \label{Opt_ARK3}\\
 \tau_k(\mA, \vb, \mAt, \vbt) = 0 \; \; \; \mbox{for} \; \; \;  k=1, . . ., P, \label{ARKtau_conditions}
 \end{align}
 \end{subequations}
 where the first three inequalities are understood component-wise, and
  $\tau_k(\mA, \vb, \mAt, \vbt)$ are the order conditions, described in the sections below.
This optimization gives the Butcher coefficients $\mA$, $\vb$, $\mAt$, and $\vbt$, and an optimal value
of the SSP coefficient  $\sspcoef = r$. Of course, for each value of $K$, a different method will be optimal.

Note that although this process defines a sufficient condition for the resulting method to be SSP with SSP coefficient  $\sspcoef = r$,
there is no reason to expect this condition to be necessary. In particular, this formulation ignores the possible interactions between $F$ and $G$
that would result in more relaxed SSP conditions.

Once again, the order conditions \eqref{ARKtau_conditions} are a key piece of the optimization, and serve as equality constraints. 
These conditions will be described in the next two subsections.

\subsubsection{Linear Order Conditions}
\noindent{\bf First order $p_{lin}=1$:} The order conditions for an additive method \eqref{eq:SC_ark_butcher_form} to be first order are 
\begin{equation}
\vb^T \ve = 1,  \; \; \; \; \; \vbt^T \ve = 1  .
\end{equation}
This has two conditions: one for the explicit part and one for the implicit part. To generalize these order conditions, we define
 $\Phi_1 = \{ \phi_{1,1}, \;  \phi_{1,2}\}=  \{ \vb^T ,  \; \vbt^T \} ,$ and the order conditions become
$\vb^T \ve = 1,  \; \; \;  \vbt^T \ve = 1  .$

\noindent{\bf Second order $p_{lin}=2$:}  Define 
$\Phi_2 = \Phi_1 \otimes \{ \mA , \; \mAt \}   = \{ \phi_{1,1} \mA, \;  \phi_{1,1} \mAt, \; \phi_{1,2} \mA, \;  \phi_{1,2} \mAt \},  $
and the order conditions are each of these right-multiplied by $\ve$ and set equal to $\frac{1}{2}$.   Recall that $\vc = \mA \ve$ and
$\vct = \mAt \ve$ 
The order conditions then become
\[ \vb^T \vc = \frac{1}{2}, \; \; \; \vb^T \vct  = \frac{1}{2}, \; \; \;
\vbt^T \vc = \frac{1}{2} , \; \; \; \vbt^T \vct = \frac{1}{2} .\]

\noindent{\bf Third order $p_{lin}=3$:}  
use $\Phi_3 = \Phi_2 \otimes \{ \mA , \; \mAt \} $ and right-multiply by $\ve$ to obtain:

\[ \vb^T \mA \vc = \frac{1}{6}, \; \; \; \vbt^T \mA \vct  = \frac{1}{6}, \; \; \; \vbt^T \mA \vc = \frac{1}{6} , \; \; \; \vb^T \mA \vct = \frac{1}{6} ,\]
\[ \vb^T \mAt \vc = \frac{1}{6}, \; \; \; \vbt^T \mAt \vct  = \frac{1}{6}, \; \; \; \vbt^T \mAt \vc = \frac{1}{6} , \; \; \; \vb^T \mAt \vct = \frac{1}{6} .\]

To obtain higher linear order, we define 
\[\Phi_{q} = \Phi_{q-1} \otimes \{\mA , \; \mAt\} =  \{ \phi_{q,1} \mA,  \phi_{q,1} \mAt, \phi_{q,2} \mA,  \phi_{q,2} \mAt, ... \phi_{q,2^q} \mA,  \phi_{q,2^q} \mAt \} \] 
and the resulting new order conditions for this order take the form
\[\phi_{q,j} \; \ve = \frac{1}{q!} \ve  \; \; \; \; \forall q, j .\]
To go from  linear order $p_{lin} = q-1$ to linear order $p_{lin}=q$ we require an additional $2^q$ order conditions. We observe that  unlike in
the implicit case, in the IMEX case the number of order conditions grows very rapidly even when both $F$ and $G$ are linear.
This is due to the many {\em coupling} order  conditions that are required for the two components to interact properly.

\subsubsection{Nonlinear Order Conditions}
%

The order conditions for first and second order are the same for both linear and nonlinear problems. 
In this section, we look at the order conditions for nonlinear orders $p^e =3,4$ for the explicit part  and $p^i=3,4$
for the implicit part.

\noindent{\bf Nonlinear order $p=3$:} For  third order, we have the linear order conditions for the explicit and implicit 
parts, and the nonlinear coupling conditions given in the section above. 
In addition, we have the  following nonlinear condition on the explicit part:
\begin{subequations}
\begin{equation}
\vb^T \vc^2 = \frac{1}{3}, \label{eq:p3_ex_nl}
\end{equation}
\mbox{the nonlinear conditions on the implicit part,}
\begin{equation}
\vbt^T \vct^2 = \frac{1}{3} \label{eq:p3_im_nl} 
\end{equation}
\mbox{and the nonlinear coupling conditions}
\begin{equation} \label{eq:p3_coupled_nla}
\vb^T \mC  \vct = \frac{1}{3}, 
\quad  \vb^T \mCt \vct = \frac{1}{3} 
\end{equation}
\begin{equation} \label{eq:p3_coupled_nlb}
 \vbt^T \mC \vc = \frac{1}{3}, 
\quad \vbt^T \mCt \vc = \frac{1}{3}.
\end{equation}
\end{subequations}



When we wish to have $p^e=3$ but use a linear $G$ (i.e $p^i = 2$) we need to include all the linear order conditions for $p_{lin} =3$, as well as
the nonlinear equations for the explicit part   \eqref{eq:p3_ex_nl}, and we must include the coupled order conditions 
in  \eqref{eq:p3_coupled_nla}. We can  neglect  the nonlinear implicit conditions  \eqref{eq:p3_im_nl} and  \eqref{eq:p3_coupled_nlb}.

\noindent{\bf Nonlinear order $p=4$:} We include all the linear and nonlinear order conditions above, for the explicit and implicit parts as
well as the coupled conditions. In addition, we have  the nonlinear order condition Equation~\eqref{eq:p4_ex_nl} for the explicit method:

\begin{subequations}
\begin{equation}
\vb^T \mA \vc^2 = \frac{1}{12} , \quad \vb^T \vc^3 = \frac{1}{4}, \quad \vb^T \mC \mA \vc = \frac{1}{8}.  \label{eq:p4_ex_nl} 
\end{equation}
\mbox{The  nonlinear order condition for the implicit part is}
\begin{equation}
\vbt^T \mAt \vct^2 = \frac{1}{12} , \quad \vbt^T \vct^3 = \frac{1}{4}, \quad \vbt^T \mCt \mAt \vct = \frac{1}{8}. \label{eq:p4_im_nl} 
\end{equation}
\end{subequations}
And the nonlinear coupled order conditions
{\small 
\begin{eqnarray} \label{eq:p4_coupled_nla}
\begin{aligned}
&  \vb^T \mA\mC\vct = \frac{1}{12} , 
\quad \vb^T \mA\mCt\vct = \frac{1}{12} ,
\quad \vb^T\mC\mCt\vc = \frac{1}{4},   
\quad  \vb^T\mC\mCt\vct = \frac{1}{4}, 
\quad \vb^T \mCt\mCt\vct = \frac{1}{4}, 
\quad \vb^T \mC \mAt \vc= \frac{1}{8}, \\
&
 \vb^T \mC \mAt \vct= \frac{1}{8},
 \quad \vb^T \mC \mA \vct= \frac{1}{8}, 
\quad \vb^T \mCt \mA \vc = \frac{1}{8}, 
\quad  \vb^T \mCt \mA \vct = \frac{1}{8} ,  
\quad \vb^T \mCt \mAt \vc = \frac{1}{8}, 
\quad \vb^T \mCt \mAt \vct = \frac{1}{8}, \\
\end{aligned}
\end{eqnarray}

\begin{eqnarray} \label{eq:p4_coupled_nlb}
\begin{aligned}
&  \vb^T \mAt\mC\vc = \frac{1}{12}, 
 \quad   \vb^T \mAt\mC\vct = \frac{1}{12}, 
  \quad \vb^T \mAt\mCt\vct = \frac{1}{12}, 
 \quad \vbt^T \mA\mCt\vc = \frac{1}{12},  
\quad \vbt^T \mA\mCt\vct = \frac{1}{12}, 
\quad \vbt^T \mAt\mC\vc = \frac{1}{12} , \\
& \vbt^T \mAt\mC\vct = \frac{1}{12},
 \quad \vbt^T \mA\mC\vc = \frac{1}{12}, 
 \quad \vbt^T\mC\mCt\vc = \frac{1}{4},   
 \quad \vbt^T\mC\mCt\vct = \frac{1}{4}, 
\quad  \vbt^T \mC \mC \vc  = \frac{1}{4},
 \quad    \vbt^T \mC \mAt \vc= \frac{1}{8},  \\
& \vbt^T \mC \mAt \vct= \frac{1}{8},
 \quad \vbt^T \mC \mA \vct= \frac{1}{8}, 
\quad \vbt^T \mCt \mA \vc = \frac{1}{8},  
\quad  \vbt^T \mCt \mA \vct = \frac{1}{8} ,  
\quad \vbt^T \mCt \mAt \vc = \frac{1}{8}, 
\quad \vbt^T \mC \mA \vc = \frac{1}{8},  \\
\end{aligned}
\end{eqnarray}
}

In the case where $G$ is linear and we wish to have $p^i=2$, we must satisfy all the linear order conditions,
the nonlinear explicit conditions \eqref{eq:p4_ex_nl}, and the coupling conditions \eqref{eq:p4_coupled_nla},
but we can neglect the nonlinear implicit conditions \eqref{eq:p4_im_nl} and the nonlinear coupling conditions
\eqref{eq:p4_coupled_nlb}

\subsection{Optimal SSP IMEX Methods} \label{IMEX_LNL_methods}
We formulate the optimization problem above in a MATLAB routine following  \cite{ketchcodes} 
and use it to generate optimized methods for different choices of $K$. 
We focus primarily on cases in which the number of stages $s$ is the same as the linear order $p_{lin}$ or
a little larger. The  nonlinear orders of interest are $p=2, 3, 4$. We note that in finding optimal methods,
we frequently observed convergence to higher order than we required or expected. For example, in the case where
we only require $p^i=2$ and $p^e>2$ the optimal SSP methods generally have $\vb=\vbt$ and $\vc=\vct$, as we would
expect from \cite{Higueras2009} 
and due to this the nonlinear implicit order is higher than required, $p^i=3$. However, we also find that in some
cases we required only $p^i=2$ and $p^e=4$, but the methods converged to $p^i=p^e=4$. For this reason,
in the following sections we describe mostly optimal methods that have $p=p^i=p^e$ (though some exceptions 
exist as is noted below). In some cases, especially where the number of stages $s$ is large and where there are more order constraints,
the optimization routine had difficulties converging and we are not comfortable stating that the methods found are optimal.
For this reason, we refer to the methods as {\em optimized} rather than optimal.

\noindent{\bf IMEX pairs for $K=\infty$}
First, we look at the case where there is no constraint resulting from the $G$ component, which is equivalent to 
setting $K=\infty$. We reformulate the optimization problem to account for this by shutting off the $\mSt$ terms 
in \eqref{Opt_ARK1} -- \eqref{Opt_ARK3}. The case of $K=\infty$ is equivalent to the case where the implicit method
is non-SSP while the explicit method is SSP. In many such cases,  researchers have focused on finding methods 
where the explicit component is SSP and the implicit component is A-stable, L-stable, or stiffly stable. These methods
were mostly of order $p\leq 3$ with few $p=4$ methods. It is known from the order barrier on explicit SSP Runge--Kutta methods
that methods of order $p>4$ cannot exist.  In our case, we focus primarily on methods that have $p_{lin} \geq 4$.
 
We first study methods with $p_{lin} \geq p$, and we do not optimize for linear stability: our rationale for studying 
 these methods is to demonstrate that where there is no SSP constraint on $G$ 
 we are able to find implicit pairs for optimal explicit SSP methods for $F$. 
 Using our optimization routine we found methods with $s$ stages, linear order $p_{lin}$ and nonlinear orders $p=p^i = p^e = 2,3$ with 
 $K=\infty$ with $s=p_{lin}, p_{lin} +1 , p_{lin}+2$.  Initially, we imposed a non-negativity conditions on all the coeffiicents, and were 
 able to find many good methods. Later, we relaxed this condition and allowed negative coefficients in $\mSt$, which allowed us to find
 additional methods.  Table \ref{IMEXKinf} contains the SSP coefficients of these methods, 
which clearly match those of the explicit LNL SSP Runge--Kutta methods found in \cite{LNL} and repeated in Table \ref{table:LNL}. 
The SSP coefficients of methods that have negative values in $\mSt$ are listed in bold.
These methods do not have large linear stability regions, as they were not optimized for this purpose. However, these
results show that the addition of an implicit component that does not have its own SSP constraint does not have an adverse effect
on the size of the SSP IMEX  method despite the addition of many new order conditions (coupling conditions) which must be satisfied.

\begin{table}
\begin{center}
\begin{tabular}{|ll|ccccc|}
\hline
nonlinear order & stages &  $p_{lin}= 3$  & $p_{lin}= 4$  & $p_{lin}= 5$  & $p_{lin}= 6$ & $p_{lin}= 7$ \\ \hline
$p=2,3$ & $s=p_{lin}$ &1.0000 &1.0000 &1.0000 & 1.0000 & 1.0000 \\ 
$p=2,3$ &$s=p_{lin}+1$  &2.0000 &2.0000 &2.0000 & 2.0000 & 2.000 \\ 
$p=2,3$ & $s=p_{lin}+2$  &2.6505 &2.6505&2.6505 & 2.6505 & {\bf 2.6505}  \\\hline
$p=4 $   & $s=p_{lin}$ & --- & --- & {\bf 0.7603} & 0.8677 & 1.0000 \\ 
$p=4$    &$s=p_{lin}+1$ & --- & 1.5082 & 1.8091   & 1.8269 & 1.9293 \\
$p=4$    & $s=p_{lin}+2$  & --- & 2.2945 & 2.5753 & 2.5629 & {\bf 2.6192}    \\ \hline
\end{tabular}
\caption{ the SSP coefficient of optimal SSP IMEX methods for nonlinear order $p= p^e=p^i=2,3,4$ with $K=\infty$. 
These match with the SSP coefficients of the explicit methods previously reported in \cite{SSPbook2011,LNL}, 
which verifies that the coupling conditions do not adversely affect the SSP coefficient in the case $K=\infty$.
Note that the  SSP coefficients in bold correspond to methods that have negative values in $\mSt$, while the other methods
have only nonnegative coefficient matrices.
 }
\label{IMEXKinf}
\end{center}
\end{table}

Next, we look for methods that  pair with popular SSP explicit Runge--Kutta methods and have large linear stability regions.
We found  a family of methods that pair with the three-stage third order Shu-Osher method:
\[ \mA=\left( \begin{array}{lll}
0 & 0 & 0 \\
1 & 0 & 0 \\
\frac{1}{4} & \frac{1}{4} & 0 \\
\end{array} \right) \qquad \vb = \left(\frac{1}{6}, \frac{1}{6}, \frac{2}{3}\right)^T.\]

These implicit pair methods are defined by $\vb =\vbt$, and 
\[ \mAt=\left( \begin{array}{rrr}
0 & 0 & 0 \\
4 \gamma+2 \beta & 1- 4 \gamma - 2 \beta & 0 \\
\frac{1}{2} - \beta - \gamma & \gamma & \beta \\
\end{array} \right)  \]
with 
\[\gamma = \frac{2 \beta^2 - \frac{3}{2} \beta  + \frac{1}{3}}{2- 4\beta} .\]
The methods are third order and for $\beta > \frac{1}{2}$ are all A-stable.
In particular, one member of this family is the nice looking method
\[ \mAt=\left( \begin{array}{rrr}
0 & 0 & 0 \\
0 & 1 & 0 \\
\frac{1}{6} & - \frac{1}{3} & \frac{2}{3} \\
\end{array} \right)  \qquad \vbt = \left(\frac{1}{6}, \frac{1}{6}, \frac{2}{3}\right)^T.\]

 Another appealing method in this family occurs  for the value $\beta= \frac{\sqrt{3}}{6} + \frac{1}{2} $.
 In this case the two non-zero values on the diagonal are equal, and so we have a type of SDIRK method
 (though with the first diagonal equal to zero).
 Note that this value of the diagonal may look familiar: in fact, it is the same as
 the diagonal value of the $s=2$, $p=3$ SDIRK method value in \cite[Table 7.2]{Hairer1},
 but that method cannot be paired with the Shu-Osher method.
 
This family of methods can be compared with two of Pareschi and Russo's methods in \cite{PareschiRusso2005}
that also pair with the  explicit SSPRK(3,3), Shu Osher method. 
The first is a three stage method that is singly implicit and L-stable and has both nonzero elements 
on the diagonal of $\tilde{A}$ but is only of overall order $p=2$. The second method is a four stage
method that is singly implicit, has four nonzero elements on the diagonal of $\tilde{A}$ but is of overall order
$p=3$. Depending on the cost of inverting the operator and the need for L-stability, the family of methods
we present, which require only two inverse solves (and possibly with the same diagonal value)
and are A-stable but not L-stable, may be preferable.

\begin{figure}[htb]
    \centering
    \begin{minipage}{0.45\textwidth}
        \centering
        \includegraphics[width=0.95\linewidth]{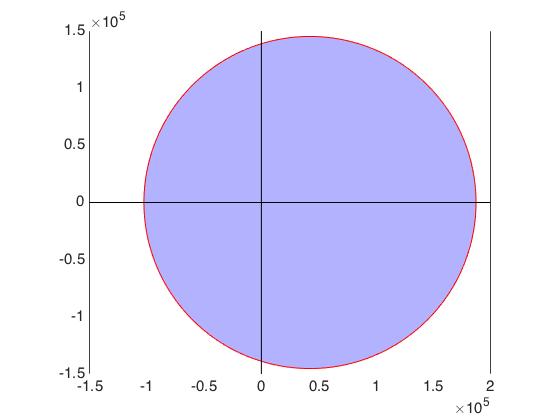} \vspace{.125in}
       \caption{Linear stability region in the complex plane of the method that pairs with Ketcheson's SSPRK(10,4) method
       which has $p^e=p_{lin}= 4$, $p^i=3$.}  
\label{Ketch104stability}
    \end{minipage}%
    \hspace{.2in}
    \begin{minipage}{0.45\textwidth}
        \centering
       \includegraphics[width=0.995\linewidth]{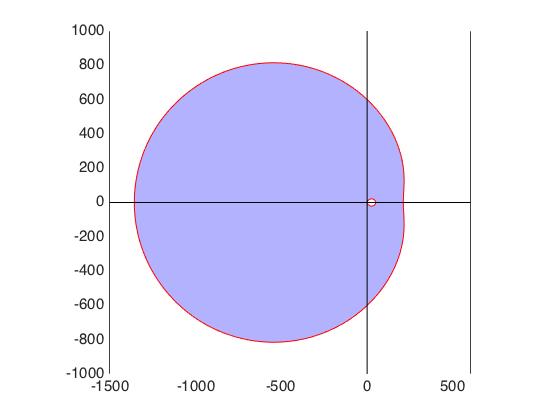}
       \caption{Linear stability region in the complex plane of the method that pairs with the LNL SSPRK with $s=10$ stages, and $p^e=4$, $p^i=3$ and $p_{lin} =6$.}
             \label{10s6pIMEX_linstab} 
           \end{minipage}
\end{figure}

A useful and efficient method is Ketcheson's  explicit SSP Runge--Kutta method with ten stages and nonlinear order $p^e=4$.
We found an SDIRK  method which pairs well with Ketcheson's SSPRK(10,4) method and has a very large linear stability region.  
The pair has linear order $p_{lin}=4$, and although it was produced by the optimizer under the assumption that $G$ is linear, we 
obtained better-than-expected nonlinear order $p^i=3$ due to the equality $\vc=\vct$.
The coefficients of this method are given in the appendix and the linear stability region is plotted in Figure \ref{Ketch104stability}. 
Note that the stability region this method crosses the real axis at $r_1 \approx -102,775$ and the imaginary axis at $r_2\approx \pm138,891$.
The optimization code was quickly able to match any value of $r= \min\{|r_1|, |r_2|\}$ that we requested. This leads us to suspect
that there may well be an A-stable method that pairs with Ketcheson's explicit SSPRK(10,4) method; however, the number of coefficients
here is prohibitive and unlike the case of the Shu-Osher SSPRK(3,3) method solving for this analytically seems difficult.

It is important to note that the explicit methods in the pair above all have an optimal SSP coefficient among methods 
of its number of stages and order, and  the coupling with an implicit method under the assumption $K = \infty$ 
does not adversely affect the SSP coefficient, so these methods are the best possible in their class
in terms of SSP coefficient.

Finally, we found an IMEX pair with $s=10$ stages, linear order $p_{lin}=6$, and nonlinear order
$p=3$. The explicit part has SSP coefficient $\sspcoef=3.3733$.
This method has an SDIRK implicit part, with a good  linear stability region as shown in  in Figure \ref{10s6pIMEX_linstab}. 
The stability region crosses the real axis at  $r_1 \approx -1,350 $ and the imaginary axis at $r_2\approx \pm 600$.

\noindent{\bf IMEX pairs for small $K$}
An interesting  case that has not been considered extensively in the literature (except indirectly in \cite{Higueras2006}),
 is the situation where the implicit component $G$ introduces a very tight time-step restriction for the strong stability 
 property to be satisfied. This is the case where the value of the parameter $K$ is very small.
We expect that in these cases the SSP coefficient of the IMEX pair will be correspondingly limited, and indeed the
 coefficients are small. Table \ref{SSPIMEX_Ktable1}
  gives the SSP coefficients of some methods we found using our optimization 
 routine.  We see in these tables that for IMEX SSP Runge--Kutta methods the SSP coefficient is reduced as
  both the linear and nonlinear orders rise; however, in general the SSP coefficients for $p=2$ and for $p=3$ are very similar,
  and for these cases the increasing linear order causes the decrease in the SSP coefficient. 
  In addition, we see that there may be some modest benefit to tailoring the method to
  the actual  value of $K$ in the problem; whether this is true in practice is investigated in  Section \ref{IMEX_SSP_verify}.
   The coefficients $\mA, \mAt, \vb$ and $\vbt$ of these methods can be downloaded at \cite{SSPIMEX_github}. 
   
 {
 \begin{table}
\begin{center}
\begin{tabular}{|c|c|ccccc|} \hline
\multicolumn{7}{|c|}{$s=p_{lin}$} \\ \hline
$ \frac{1}{K}$ & p  & \multicolumn{5}{|l|}{$p_{lin}$}  \\  
                      &     & 3 & 4 & 5 & 6 & 7 \\ \hline
10   	& 2 	 &  2.030$\times 10^{-1}$ & 1.729$\times 10^{-1}$& 1.520$\times 10^{-1}$&  1.323$\times 10^{-1}$ &  1.109$\times 10^{-1}$ \\
	& 3   &   1.492$\times 10^{-1}$ & 1.727$\times 10^{-1}$ & 1.520 $\times 10^{-1}$ &1.323$\times 10^{-1}$ & 0.932$\times 10^{-1}$ \\
	& 4    &  -- & -- & --		    	&  1.132$\times 10^{-1}$ & 0.849$\times 10^{-1}$ \\ \hline
100 	& 2	 &  2.23$\times 10^{-2}$ & 2.06$\times 10^{-2}$ &  1.58$\times 10^{-2}$ &  1.51$\times 10^{-2}$ &  0.99$\times 10^{-2}$ \\
	& 3 	&   1.63$\times 10^{-2}$ & 2.05$\times 10^{-2}$ &  1.58$\times 10^{-2}$ &  1.51$\times 10^{-2}$ &  0.99$\times 10^{-2}$ \\
	& 4    &  -- & -- & --& 1.18$\times 10^{-2}$ &  0.89$\times 10^{-2}$   \\ \hline 
\end{tabular} 
 \begin{tabular}{|c|c|ccccc|} \hline
\multicolumn{7}{|c|}{$s=p_{lin}+1$} \\ \hline
$ \frac{1}{K}$ & p  & \multicolumn{5}{|l|}{$p_{lin}$}  \\ 
                      &     & 3 & 4 & 5 & 6 & 7 \\ \hline
10   	& 2 	 &  3.570$\times 10^{-1}$ & 3.099$\times 10^{-1}$ & 2.891$\times 10^{-1}$ & 2.307$\times 10^{-1}$ & 1.756$\times 10^{-1}$ \\
	& 3   &   2.837$\times 10^{-1}$ & 3.084$\times 10^{-1}$ & 2.767$\times 10^{-1}$ & 2.288$\times 10^{-1}$ & 1.756$\times 10^{-1}$ \\
	& 4    &-- 		& 2.001$\times 10^{-1}$ & 2.164$\times 10^{-1}$ & 1.986$\times 10^{-1}$ & 1.755$\times 10^{-1}$  \\ \hline
100 	& 2	&  3.97$\times 10^{-2}$ & 3.65$\times 10^{-2}$ & 3.26$\times 10^{-2}$ & 2.98$\times 10^{-2}$ & 1.85$\times 10^{-2}$ \\
	& 3 	&  3.23$\times 10^{-2}$ & 3.56$\times 10^{-2}$ & 3.20$\times 10^{-2}$ & 2.98$\times 10^{-2}$ & 1.79$\times 10^{-2}$  \\
	& 4    &  -- 		& 2.26$\times 10^{-2}$ & 2.28$\times 10^{-2}$ & 1.94$\times 10^{-2}$ & 1.66$\times 10^{-2}$  \\ \hline
\end{tabular}
\end{center}
\caption{SSP coefficients of some optimized IMEX methods with $s=p_{lin}$ and  $s=p_{lin}+1$
 stages optimized for different values of $K$.}
\label{SSPIMEX_Ktable1}
\end{table}
}
 
%
%

      These methods are not at all optimized in terms of the linear stability region, as the main constraint of interest for these
 methods comes from the  SSP condition, and indeed their linear stability region is similar to explicit methods. 
 We believe these methods may be useful where the cost of the implicit solution is negligible
 while the SSP condition is needed for small values of $K$. However, even in cases where these methods are not directly useful,
 they give us an upper bound on the possible SSP coefficients of SSP IMEX methods for typical values of $K$, and so serve as a guide
 to what we can expect when we co-optimize these methods with other desirable properties.  
 The numerical tests, particularly those in  Section \ref{IMEX_SSP_verify}, show the performance of some of these
 methods for the types of problem that require the SSP property for both the explicit $F$ and implicit $G$.

\subsection{Numerical Experiments} \label{Numerical2}

\subsubsection{Convergence studies} \label{convergence} 
To test the accuracy of the SSP IMEX Runge--Kutta methods, we consider a linear or nonlinear explicit part. As our main motivation for these
methods are where the implicit component is linear, we test the methods on problems where $G$ is linear.
 In these tests we confirm that our methods give us the expected linear and nonlinear orders
for the explicit and implicit parts.

\begin{figure}[htb]
\includegraphics[width=0.475\linewidth]{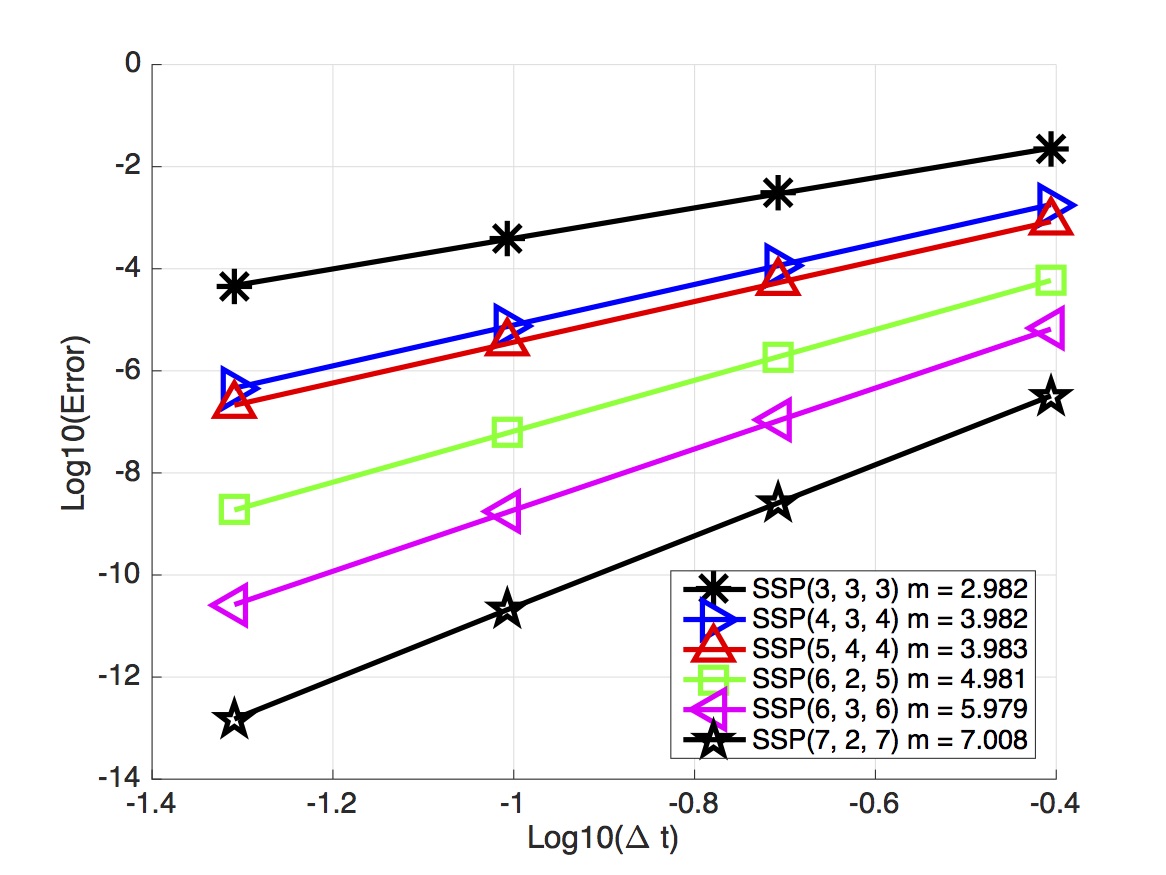}
\includegraphics[width=0.475\linewidth]{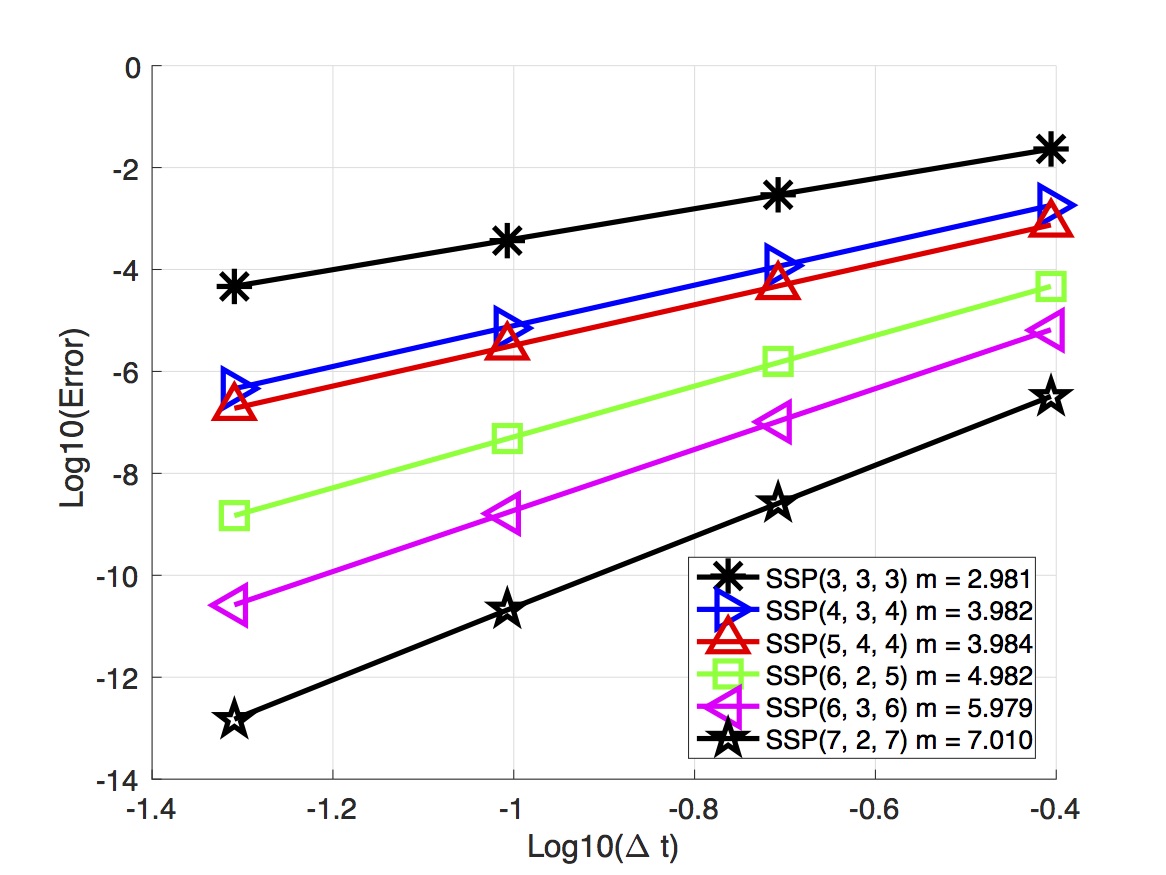}
       \caption{Convergence plots for an advection diffusion problem in Example 3.1. 
       Methods shown are $(s, p^i,p^e,p_{lin})  =(3,3,3), (4,3,4), (5,4,4), (6,2,5),(6,3,6), (7,2,7)$.
       On the left are the methods created for $K=\infty$ and on the right for $K=10^{-2}$.}
            \label{fig:IMEXconv_lin_lin}
\end{figure}

 \noindent{\bf Example 3.1: convergence study with explicit linear advection and implicit linear diffusion.}
We approximate the solution to the equation \[ U_t +  U_x = \epsilon U_{xx}  \; \; \; \; \; \;  x \in [0, 2 \pi] \]
where $\epsilon = 0.01$ with periodic boundary conditions and sine wave initial conditions. 
 We discretize the spatial grid with $N=8$ equidistant points and use the Fourier pseudospectral 
differentiation matrix  for both the first and second derivative \cite{HGG2007}.
The advection part is dealt with explicitly and the diffusion implicitly (and exactly using MATLAB's backslash operator).
 Temporal grid refinement with pseudospectral approximation of the spatial derivative.
We use a range of time steps, $\Delta t = \lambda \Delta x$, where $\Delta x = \frac{\pi}{4}$ and 
 we pick $\lambda = \frac{1}{16},\frac{1}{8}, \frac{1}{4}, \frac{1}{2}$ to compute the solution to final time $T_f = 5.0$.
The methods we test are our IMEX SSP methods with $s$ stages, nonlinear order $p=p^i=p^e$ and linear order $p_{lin}$
$(s, p,p_{lin} ) =(3,3,3), (4,3,4), (5,4,4), (6,2,5),(6,3,6), (7,2,7)$ designed for $K=\infty$ and for $K=10^{-2}$.
The $\ell_2$ errors are measured compared to the exact solution and shown in 
Figure \ref{fig:IMEXconv_lin_lin}. We note that the methods designed for $K=\infty$ and those designed for $K=10^{-2}$ 
perform essentially  the same on convergence studies. The orders (slope $m$ of the line) are taken by taking a linear fit 
using MATLAB's {\tt polyfit} and given in the caption.
We see that, as expected, the linear design order is apparent in this linear problem.  
It is interesting to note that the $(s,p,p_{lin} ) =(3,4,4)$ and the $(s, p,p_{lin})  =(4,4,5)$ have the same slope, because their linear order
is the same, but the latter method has a smaller error constant so the errors are smaller.

\begin{figure}[htb]
\includegraphics[width=0.475\linewidth]{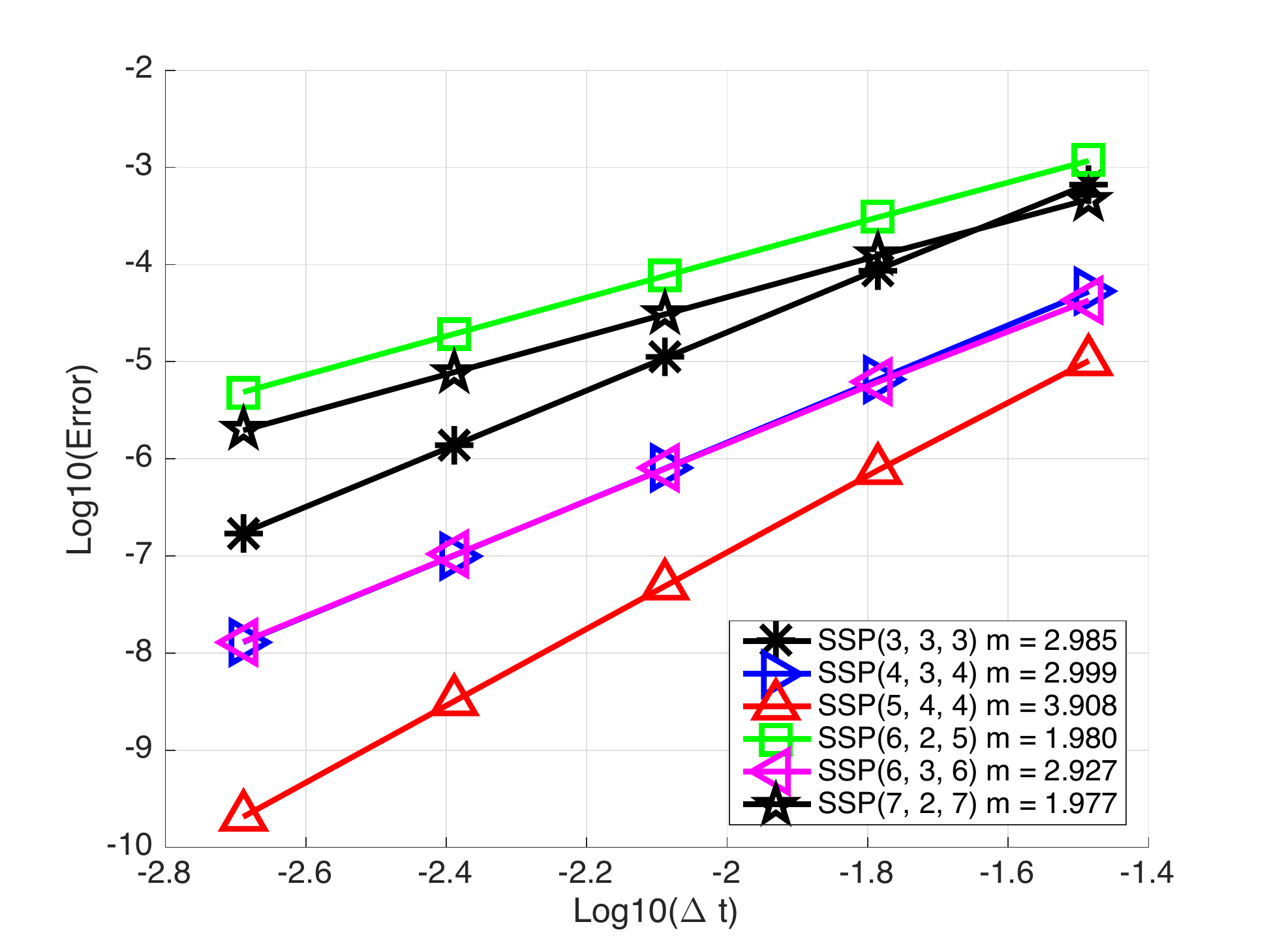}
\includegraphics[width=0.475\linewidth]{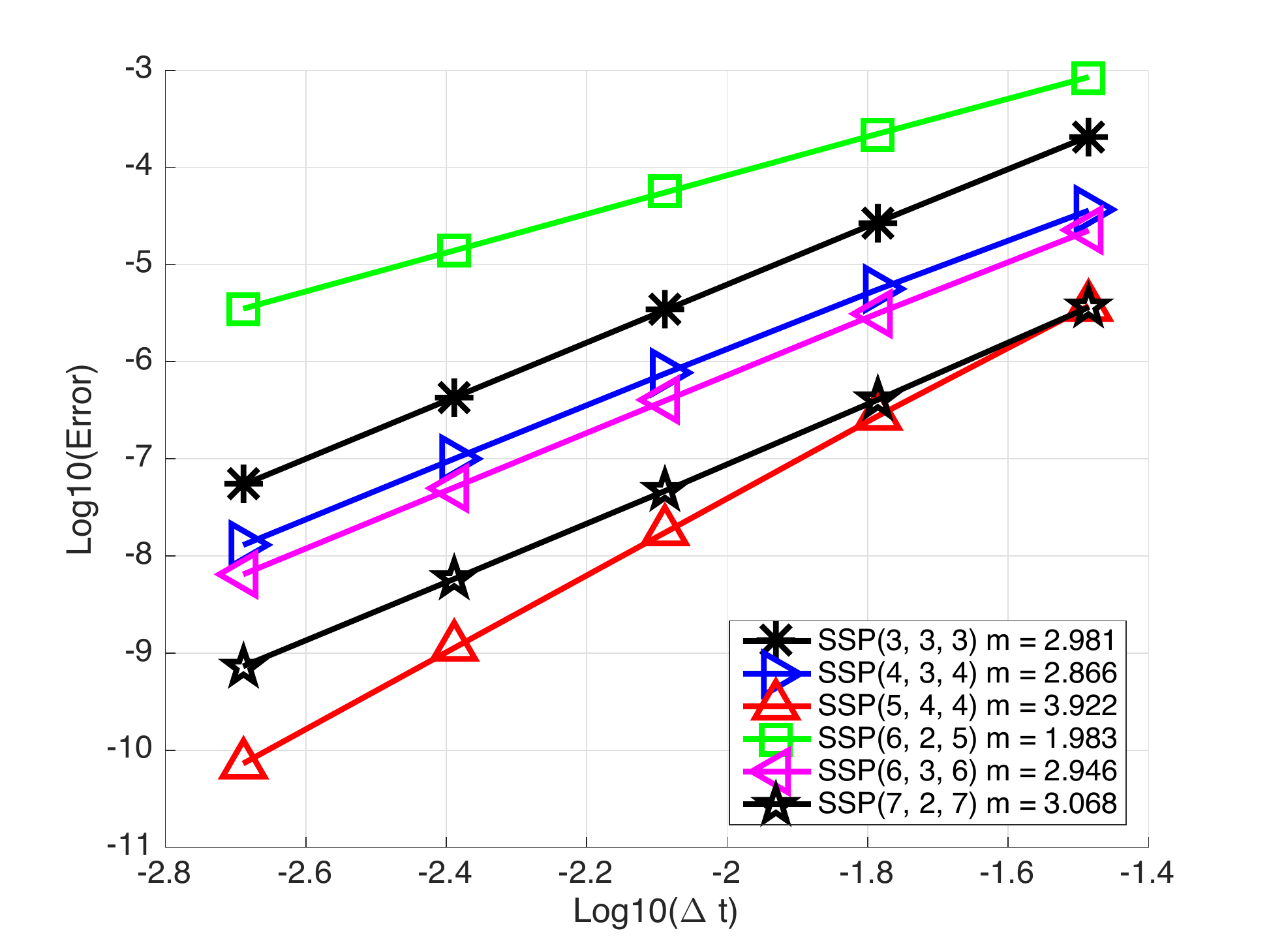}
       \caption{Convergence plots for an the explicit Burgers' with implicit linear advection 
       problem in Example 3.2. 
       Methods shown are $(s, p,p_{lin} ) =(3,3,3), (4,3,4), (5,4,4), (6,2,5),(6,3,6), (7,2,7)$.
        On the left are the methods created for $K=\infty$ and on the right for $K=10^{-2}$.
       }
            \label{fig:IMEXconv_nl_lin}
\end{figure}

\noindent{\bf Example 3.2: convergence study with explicit Burgers' and implicit linear advection.}
We approximate the solution to the equation \[ U_t +  \left( \frac{1}{2} U^2 \right)_x + U_x = 0  \; \; \; \; \; \;  x \in [0, 2 \pi] \]
with periodic boundary conditions and sine wave initial conditions. 
We discretize the spatial grid with $N=24$ equidistant points and use the Fourier pseudospectral 
differentiation matrix  \cite{HGG2007}.
The Burgers' flux $\left( \frac{1}{2} U^2 \right)_x$
is dealt with explicitly and the linear advection implicitly (and exactly using MATLAB's backslash operator).
We perform  grid refinement with pseudospectral approximation of the spatial derivative,
we use a range of time steps, $\Delta t = \lambda \Delta x$
where we pick $\lambda = \frac{1}{128}, \frac{1}{64}, \frac{1}{32} , \frac{1}{16}, \frac{1}{8} $ to 
compute the solution to final time $T_f = 0.8$ before the shock forms.
The methods we test are our IMEX SSP methods with
$(s, p,p_{lin} ) =(3,3,3), (4,3,4), (5,4,4), (6,2,5),(6,3,6), (7,2,7)$ 
designed for $K=\infty$ and for $K=10^{-2}$.
The $\ell_2$ errors are measured compared to the the results from MATLAB's {\tt ode45} and graphed  in 
Figure \ref{fig:IMEXconv_nl_lin}. The orders (slope $m$ of the line) are taken by taking a linear fit using MATLAB's {\tt polyfit}
and given in the caption. We see that the nonlinear order $p$ typically dominates: this occurs since the
linear order $p_{lin}$ is no smaller than the nonlinear orders. In the case $(7,2,7)$ for $K=10^{-2}$
we get  a better order than expected by one, due to a small error constant  ($\approx 10^{-8}$) for this method.

%

\subsubsection{Numerical verification of the SSP properties of these methods} \label{IMEX_SSP_verify}

In this section we focus on the behavior of our time-stepping methods on problems which require the 
SSP property. In particular, we use out IMEX SSP methods for problems where the desired  property 
to be preserved is the total variation diminishing (TVD) property and where the spatial discretization is simple
enough that the theoretical bound on the TVD condition is easy to see.

\noindent{\bf Example 4.1: Explicit linear advection with  implicit linear advection. }
Consider  the linear equation 
\[ U_t +  U_x + 100 U_x = 0  \; \; \; \; \; \;  x \in [-1, 1] \]
with periodic boundary conditions and step function initial conditions
 \begin{equation}
\label{eqn:sq_wave_ic}
	u_0(x) = \left\{ \begin{array}{ll}
		1 & \text{if}\ \frac{1}{4} \leq x \leq \frac{1}{2}, \\
		0 & \text{otherwise},
\end{array} \right. 
\end{equation}
We use a first order upwind differencing the  linear advection terms,
with  $N=301$ points in space. The time-step is set to $\dt = \lambda \dx$, where different values of 
$\lambda$ are tested to observe the value at which the discretization is no longer TVD
 (i.e. the total variation rises by more than $10^{-10}$ at any given time-step).  We first compared the
 Predicted $\dt$ and  Observed $\dt$  for which the TVD property is
violated. The results are very similar to those in Table \ref{NonlinearDT} and so we do not report them here.

We consider several IMEX methods that were optimized for $K=1, \frac{1}{10}, \frac{1}{100}$ 
and see how their observed time-step for the TVD to be satisfied compares to the time-step predicted 
by the theory, even when the methods were designed for different values of $K$. 
We also compare the performance of these methods to the situation where both components are 
advected by the optimal explicit LNL SSP Runge--Kutta method of the corresponding number of 
stages, nonlinear order, and linear order. For completeness, we also include a comparison to
treating both components implicitly using an optimal implicit LNL SSP Runge--Kutta method.

Table \ref{LinearDT}  contains the observed and predicted time-step for which the TVD property is
violated. Clearly, the observed $\dt$ is never smaller than the predicted $\dt$, as the 
SSP property  provides a guarantee of this property. We note that while the predicted and observed TVD time-step 
are always {\em exact} when both components are treated implicitly, and frequently so when they are both
treated explicitly, this is not as often the case where we use an IMEX method (though in many cases, particularly those
noted in Table  \ref{NonlinearDT}, we do have sharp agreement). It is notable that using an IMEX method 
that was optimized for a  value of $K$ close to the actual value of $\frac{1}{\omega}$ we generally see larger 
allowable time-step, though this difference may not be large enough to generate optimal methods for many different values.
The main question we wish to answer is whether these methods are useful as replacements to the fully explicit
method. In these cases we see that using an IMEX scheme allows us to take a time-step that is between $1.5$ and $2$
times the fully explicit time-step (see the ratio column). There are very limited circumstances in which this  larger time-step
is truly enough to offset the cost of the implicit solver; however, this is what one would have expected from the fact that
the SSP coefficient of the explicit methods are bounded by the number of stages while the 
SSP coefficient of the implicit methods are bounded by twice the number of stages. The more interesting 
observation  in this table may be the ratio of the fully implicit allowable time-step to the fully explicit allowable time-step:
these can be quite large, possibly large enough to offset the implicit solver cost in some cases.

\begin{table}[t] 
\hspace{0.8in} \begin{tabular}{|cccccccc|}  \hline 
 Type & p & $p_{lin}$ & s &$ \frac{1}{K}$ & Predicted $\dt$ & Observed $\dt$  & ratio  \\ \hline
Explicit  & 2  & 5 &  6  & --     & 1.98$\times 10^{-2} $ &  1.98$\times 10^{-2} $ &  --\\
IMEX  &       &    &      &   1    & 1.25$\times 10^{-2}$ &  3.07$\times 10^{-2} $ & 1.55\\
IMEX  &  &    &      &  10   &  2.89$\times 10^{-2}$ &  3.33$\times 10^{-2} $ & 1.68\\
IMEX   & &    &      &  100 &  3.26$\times 10^{-2}$    &  3.37$\times 10^{-2} $ &  1.70 \\  
Implicit & &    &      &  -- & 6.59$\times 10^{-2}$ & 6.59$\times 10^{-2}$ & 3.32 \\ \hline
Explicit  &  2   &  6  &   6     &     	--    & 9.90$\times 10^{-3} $     &	9.90$\times 10^{-3} $  &   -- \\ 
IMEX  & &    &    	&	1     &  6.00$\times 10^{-3} $  &  1.47$\times 10^{-2} $ & 1.48 \\
IMEX  & &    &      &	10   &  1.32$\times 10^{-2} $ &    1.94$\times 10^{-2} $  & 1.96 \\
IMEX  & &    &   	&	100 &  1.51$\times 10^{-2} $   &  1.97$\times 10^{-2} $  & 1.99 \\
Implicit & &    &      &  -- &  5.12$\times 10^{-2} $ & 5.12$\times 10^{-2} $ & 5.17   \\ \hline
Explicit   &  3  &  5   &  5    &      	-- &   9.90$\times 10^{-3} $ &     9.90$\times 10^{-3} $ & --  \\
IMEX  & 	   &       &        &	       1    &  6.34$\times 10^{-3} $ &    2.00$\times 10^{-2} $ & 2.02 \\
IMEX  & 	   &       &        &	      10   &   1.52$\times 10^{-2} $ &   2.24$\times 10^{-2} $ & 2.26 \\
IMEX  & 	   &       &        &	     100  &    1.58$\times 10^{-2} $  &  2.38$\times 10^{-2} $ & 2.40 \\ 
Implicit & &    &      &  -- & 4.96$\times 10^{-2} $ & 4.96$\times 10^{-2} $ &  5.01 \\  \hline
Explicit  &  3  &   5   &   6    &      	--    &     1.98$\times 10^{-2} $ &    1.98$\times 10^{-2} $    & -- \\
IMEX  & &    &    	&		1     &	1.23$\times 10^{-2} $ &   3.10$\times 10^{-2} $  & 1.56 \\
IMEX  & &    &      & 	 		10    &   2.77$\times 10^{-2} $   & 3.30$\times 10^{-2} $        & 1.66 \\
IMEX  & &    &      &			100  &  3.20$\times 10^{-2} $    &3.50$\times 10^{-2} $         & 1.76 \\ 
Implicit & &    &      &  -- & 6.59$\times 10^{-2} $ &  6.59$\times 10^{-2} $ &   3.32  \\   \hline
Explicit   &  4  &  6   &  6    &      	-- &   8.59$\times 10^{-3} $    & 9.90$\times 10^{-3} $  & --  \\
IMEX  & 	   &       &        &	   1      & 	5.12$\times 10^{-3} $    & 1.91$\times 10^{-2} $  & 1.93 \\
IMEX  & 	   &       &        &	  10      & 1.13$\times 10^{-2} $     &  2.14$\times 10^{-2} $  & 2.16 \\
IMEX  & 	   &       &        &	  100     &	1.18$\times 10^{-2} $    &  2.09$\times 10^{-2} $  &  2.11 \\ 
Implicit & &    &      &  -- & 5.09$\times 10^{-2} $ &  5.09$\times 10^{-2} $ &   5.92 \\ \hline
\end{tabular} 
\caption{Comparison of the theoretical and observed allowable time-step before an SSP violation 
of $10^{-12}$ occurs for the linear problem  in Example 4.1
with wavespeed $\omega=100$ and for method optimized for the value of $K$ above.}
\label{LinearDT}
\end{table} 

\noindent{\bf Example 4.2: Explicit Burgers' with implicit linear advection. }
We consider  the equation 
\[ U_t +  \left( \frac{1}{2} U^2 \right)_x + \frac{1}{K} U_x = 0  \; \; \; \; \; \;  x \in [-1, 1] \]
with periodic boundary conditions and step function initial conditions
 \begin{equation}
\label{eqn:sq_wave_ic}
	u_0(x) = \left\{ \begin{array}{ll}
		1 & \text{if}\ \frac{1}{4} \leq x \leq \frac{1}{2}, \\
		0 & \text{otherwise},
\end{array} \right. 
\end{equation}
We use a first order upwind differencing for both the Burgers' term and the linear advection term,
with  $N=301$ points in space. The methods used are those optimized for the  value of $K$ that corresponds to the
wavespeed $\omega=\frac{1}{K}$, and these IMEX methods have
$s$ stages, nonlinear order $p^e=p^i=p$ and nonlinear order $p_{lin}$.
The time-step is set to $\dt = \lambda \dx$, where different values of 
$\lambda$ are tested to observe the value at which the discretization is no longer TVD.  This value of $\dt$ is then called
the ``Observed $\dt$'' while the value we would expect from the SSP coefficient is called the 
``Predicted $\dt$''.

Table \ref{NonlinearDT} contains the observed and predicted time-step for which the TVD property is
violated. Clearly, the observed $\dt$ is never smaller than the predicted $\dt$, as the 
SSP property  provides a guarantee of this property. As seen in the table, in some cases, the 
methods  feature close agreement between the predicted and observed SSP coefficient, while in others
the actual time-step allowed for this problem is greater than predicted. We note that these
methods performed very similarly on the linear problem in Example 4.2 above, so the results were not repeated
in a separate table. These results, and in particular the sharpness of the predicted allowable time-step in 
several cases, demonstrate that the theoretical SSP property is a good predictor of the actual behavior of the method.

\begin{table}[t] 
\hspace{0.8in} \begin{tabular}{|cccccc|}  \hline 
 p & $p_{lin}$ & s &$ \omega=\frac{1}{K}$ & Predicted $\dt$ & Observed $\dt$ \\ \hline
 2 & 4 & 5 & 10   & 3.099$\times 10^{-1} $  & 3.099$\times 10^{-1} $   \\
    &     &   & 100 &  3.65$\times 10^{-2} $  & 3.65$\times 10^{-2} $   \\ \hline
 3 &  4  & 5 & 10  & 3.084$\times 10^{-1} $   & 3.084$\times 10^{-1} $   \\
    &     &   & 100 &  3.56$\times 10^{-2} $   & 3.56$\times 10^{-2} $   \\ \hline
 3 &  5 & 5 &  10 & 1.520$\times 10^{-1} $  &  2.135$\times 10^{-1} $ \\
    &     &   & 100 & 1.58$\times 10^{-2} $    & 2.39$\times 10^{-2} $  \\ \hline
 2 & 4 & 6 &  10  & 4.088$\times 10^{-1} $  & 4.088$\times 10^{-1} $  \\ 
     &     &   & 100 & 4.67$\times 10^{-2} $  & 4.67$\times 10^{-2} $ \\ \hline
 3  & 4  & 6 & 10 & 4.171$\times 10^{-1} $  & 4.171$\times 10^{-1} $ \\
     &     &   & 100 & 4.59$\times 10^{-2} $ &  4.59$\times 10^{-2} $ \\\hline
 3  & 6  & 7 & 10 & 2.29$\times 10^{-1} $ &    3.00$\times 10^{-1} $ \\
     &     &   & 100 & 2.98$\times 10^{-2} $ &    3.50$\times 10^{-2} $ \\ \hline
 4  &  6   &  6    &  10 & 1.13$\times 10^{-1} $ &    2.01$\times 10^{-1} $   \\
   &     &   & 	100 &    1.18$\times 10^{-2} $ &   2.10$\times 10^{-2} $    \\ \hline
\end{tabular} 
\caption{Comparison of the theoretical and observed SSP coefficients that preserve the nonlinear stability properties in Example 4.2
with the methods optimized for the value $\frac{1}{K}$ to match the wavespeed $\omega$.}
\label{NonlinearDT}
\end{table} 

\vspace{-.15in}
\section{Conclusions}
\vspace{-.1in}

In this work, we investigated implicit and IMEX SSP methods with very high linear order.
We first considered implicit SSP Runge--Kutta methods that have high order for linear problems and order $p=2$ 
for nonlinear problems.  We show that  as is the case with explicit methods we are able to find implicit methods  with $p_{lin} > 6$
with no linear order barrier. We were further able to show that (as in \cite{LNL}) the optimal methods had no better SSP coefficient 
that optimal methods that are diagonally implicit with a low storage formulation. Thus we can say that optimal implicit SSP Runge--Kutta methods
for any linear order $p_{lin}$ and nonlinear order $p=2$  are diagonally implicit and low storage.

We continued our study of implicit  linear/nonlinear (LNL) methods for larger nonlinear orders $p>2$
following the approach used in \cite{LNL}. We designed SSP implicit Runge--Kutta methods where the 
linear order exceeded the nonlinear order (called LNL methods),  and observed that  the SSP coefficient 
for methods with linear order $p_{lin} \leq 9$ was similar whether we chose  nonlinear order $p=2$ and $p=3$, 
and that if we went up to $p=4$, the SSP coefficient was not typically significantly reduced as long as we had
sufficiently many stages.  This implies that just as the case for the explicit LNL methods, as we increase the 
number of stages and maintain a moderate nonlinear order $p=3,4$, the linear order is the main constraint 
on the SSP coefficient. However, if we increase the nonlinear order to $p=5,6$, this  did have a strong adverse 
effect on the size of the SSP coefficient. We also observed that the SSP coefficients of the implicit LNL SSP Runge--Kutta methods
were up to six times as large as those of the corresponding explicit LNL SSP Runge--Kutta of the same number of stages and linear
and nonlinear orders.
We report the coefficients of  methods with high linear order
 and moderate number of stages that have reasonable SSP coefficients, for example the  LNL implicit RK 
 $(s,p,p_{lin})=(6,4,6)$ method with SSP coefficient $\sspcoef=0.856$. These methods may be useful when a 
 high linear order method is desired, while still reasonable for use with nonlinear problems.
We then verified the convergence of these methods as well as the sharpness of the SSP coefficient on sample problems.

The second part of this paper focused on implicit-explicit methods, where as above the methods have higher linear order than nonlinear order.
These methods are of value when we desire both the explicit and implicit parts to have the SSP property, and are particularly beneficial 
when the linear part of a problem is treated implicitly, and the nonlinear part explicitly.
We found diagonally implicit non-SSP Runge--Kutta methods with large regions of linear stability that pair with well-known  
explicit SSP  methods. Next, we found pairs of IMEX SSP methods, with high linear order,  that were optimized for the relationship
between the forward Euler conditions of each component.
We verified the order of these methods on a variety of problems, and the sharpness of the SSP coefficient on a typical test cases, 
and conclude that these IMEX LNL methods demonstrate high order convergence and  perform as desired in terms of SSP.

This work shows that it is possible to produce implicit and IMEX SSP methods of very high order if we consider only the linear order.
For implicit methods, we found methods of up to order $p_{lin}=9$ and for the IMEX methods up to order $p_{lin}=7$, which show that
the order barriers of the implicit and IMEX SSP Runge--Kutta methods apply only to the  nonlinear order. 
The approach described in this work can be used to produce optimal SSP implicit and IMEX methods with additional 
desired properties, as well. The SSP coefficients presented in this work give a baseline with which to compare any
methods that are SSP and also have other properties, and give a value to strive toward. 
As expected at the outset, the allowable time-step for these methods is very restricted, especially for the IMEX schemes and therefore
the cost of the implicit solvers is still the major issue that arises in the solution of these methods. 
However, in particular cases (such as where the implicit term is a constant coefficient linear term)  where the cost of the implicit solver 
is controllable,   these methods may be of more than theoretical value.  

{\bf Acknowledgements} This work was supported by AFOSR grant FA9550-15-1-0235, and 
 was partially supported by DOE NNSA ASC Algorithms effort, the DOE Office of Science AMR program 
 at Sandia National Laboratory under contract DE-AC04-94AL85000

\appendix
\section{Coefficients of some optimized  methods}
\subsection{Methods for $K=\infty$}
\subsubsection{An IMEX pair based on Ketcheson's SSPRK10,4 method}
Ketcheson's SSPRK10,4 method is given by a vector $\vb$ where $(\vb)_i = \frac{1}{10}$ for $  1 \leq i \leq 10$
and the coefficients in $\mA$ are \\
\begin{eqnarray*}
 a_{i,j}  =  \frac{1}{6}  & \mbox{for} \; \;  j < i \leq 5 \\
 a_{i,j}  =   \frac{1}{15} & \mbox{for} \; \;  1 \leq j  \leq 5  \; \; \mbox{and} \; \;  6 \leq i  \leq 10 \\
 a_{i,j}  =  \frac{1}{6}  & \mbox{for} \; \;  6 \leq j < i  \leq 10
 \end{eqnarray*}

The SDIRK method described above that pairs with Ketcheson's 10,4 method has coefficients $\vbt = \vb$ and $\vct = \vc$
and  the non-zero coefficients in $\mAt$ are: 

\smallskip

\begin{tabular}{llll}
   \multicolumn{3}{c}{ $\tilde{a}_{i,i}=0.929729066567767$  \; for \;  $  2 \leq i \leq 10$}   \\
$\tilde{a}_{2,1}=$ -0.763062399901101 & $\tilde{a}_{3,1}=$ -1.929471352156769 &  $\tilde{a}_{3,2}=$1.333075618922335   \\
$\tilde{a}_{4,1}=$ -1.746903568350466 & $\tilde{a}_{4,2}=$0.408445589167274  & $\tilde{a}_{4,3}=$0.908728912615425  \\
$\tilde{a}_{5,1}=$0.565228647234277  & $\tilde{a}_{5,2}=$ 1.133923847131481  & $\tilde{a}_{5,3}=$-1.557731112458759 \\
$\tilde{a}_{5,4}=$-0.404483781808100 & $\tilde{a}_{6,1}=$1.982844041162849  & $\tilde{a}_{6,2}=$-1.490145231639306 \\
$\tilde{a}_{6,3}=$ -0.008867539995790 & $\tilde{a}_{6,4}=$ -1.160584799688216 & $\tilde{a}_{6,5}=$0.080357796926028 \\
$\tilde{a}_{7,1}=$ 0.221597237328096 &  $\tilde{a}_{7,2}=$1.616180514391033 &$ \tilde{a}_{7,3}=$0.142461646204330 \\
$\tilde{a}_{7,4}=$-0.868274370597692 & $\tilde{a}_{7,5}=$-1.991484177541085 & $\tilde{a}_{7,6}=$0.449790083647550 \\
$\tilde{a}_{8,1}=$  -1.546919287943971 & $\tilde{a}_{8,2}=$  1.854908818861482 & $\tilde{a}_{8,3}=$1.205736394483380 \\
$\tilde{a}_{8,4}=$ -0.314106013195022  & $\tilde{a}_{8,5}=$ 0.915344917019776 & $\tilde{a}_{8,6}=$-1.386044641065531 \\
$\tilde{a}_{8,7}=$ -0.991982588061215  &  $\tilde{a}_{9,1}=$-0.091706218761790 & $\tilde{a}_{9,2}=$1.633885494435077   \\
$\tilde{a}_{9,3}=$0.932276645625014   & $\tilde{a}_{9,4}=$ -1.944938658929756  & $\tilde{a}_{9,5}=$-1.977191163021469  \\
$\tilde{a}_{9,6}=$ 1.963551314474635  & $\tilde{a}_{9,7}=$ -1.871583791474667 & $\tilde{a}_{9,8}=$1.259310644418523   \\
$\tilde{a}_{10,1}=$-1.527363916489275 & $\tilde{a}_{10,2}=$1.982728522581499  & $\tilde{a}_{10,3}=$1.859310770893058  \\
$\tilde{a}_{10,4}=$-1.881872618524453  & $\tilde{a}_{10,5}=$ 1.047237251794738 & $\tilde{a}_{10,6}=$-1.831562507581245 \\
$\tilde{a}_{10,7}=$1.992738025048269   & $\tilde{a}_{10,8}=$-1.135512580190266  &$\tilde{a}_{10,9}=$-0.435432014100091 \\
    \end{tabular}

    \subsection{Methods for $K=10$}
     \subsubsection{An IMEX pair with $p^e=p^i=3$, $s=p_{lin}=5$}
     The five stage IMEX pair with 
     
        The nonzero coefficients in $\mA$ are given by: 
        
\begin{tabular}{llll}
 $a_{2,1}=$  0.740010097277110  & $a_{3,1}=$  0.058133047039451 & $a_{3,2}=$ 0.516728366555161    \\
 $a_{4,1}=$  0.327995830636910  &  $a_{4,2}=$ 0.028076226778328  &  $a_{4,3}=$  0.357399140460949    \\
  $a_{5,1}=$ 0.255837111227683   &  $a_{5,2}=$0.074862387600713   &  $a_{5,3}=$0.116959465282915  \\
   $a_{5,4}=$ 0.195688888775226    & & \\
    \end{tabular}

\begin{tabular}{llll}
$\tilde{a}_{2,1}=$ 0.583773436668528   & $\tilde{a}_{2,2}=$0.156236660608582  & $\tilde{a}_{3,1}=$0.276599046373025\\
$\tilde{a}_{3,2}=$ 0.012273492120642   & $\tilde{a}_{3,3}=$ 0.285988875100944  & $\tilde{a}_{4,1}=$  0.348206780427965 \\
$\tilde{a}_{4,2}=$ 0.349725300350930   & $\tilde{a}_{4,3}=$0.015539117097292   & $\tilde{a}_{5,1}=$ 0.226390976173007   \\
$\tilde{a}_{5,2}=$0.140957344725959    & $\tilde{a}_{5,3}=$0.080310643212345   & $\tilde{a}_{5,4}=$0.195688888775226 \\
    \end{tabular}  
   
             \[    \mbox{and} \; \; \; \;    \vb = \vbt=   \left(   \begin{array}{l} 
      0.243859806139543\\
   0.180742612023724\\
   0.161824368384123\\
   0.101972004412874\\
   0.311601209039737\\
    \end{array} \right).  \]
   
         \subsubsection{An IMEX pair with $p^e=p^i=4$, $p_{lin}=6$}
   The nonzero coefficients in $\mA$ are given by: 
        
\begin{tabular}{llll}
$a_{2,1}=$  0.376055593238192  &   $a_{3,1}=$0.127359848364171  &  $a_{3,2}=$0.318823868640133    \\
$a_{4,1}=$ 0.184561142322538   & $a_{4,2}=$ 0.021362084173389  & $a_{4,3}=$  0.297673384659880    \\
$a_{5,1}=$ 0.132655730337691   & $a_{5,2}=$ 0.006786263788702  & $a_{5,3}=$ 0.094405536658542  \\
$a_{5,4}=$ 0.228410568816280    & $a_{6,1}=$ 0.114983915662321  & $a_{6,2}=$ 0.136226885295266  \\
$a_{6,3}=$ 0.045369437546957    & $a_{6,4}=$ 0.109769611285921   & $a_{6,5}=$0.326960155246028    \\
$a_{7,1}=$ 0.122086625034326   &  $a_{7,2}=$0.097697571022518   &  $a_{7,3}=$0.158995977454046 \\
$a_{7,4}=$ 0.117285498485044  &  $a_{7,5}=$ 0.211659248630559 &  $a_{7,6}=$  0.261454998381366  \\
    \end{tabular}

\noindent The nonzero coefficients in $\mAt$ are given by: 
 
 \begin{tabular}{llll}
$\tilde{a}_{2,1}=$ 0.376055593238191   &$\tilde{a}_{3,1}=$ 0.158832832190656   & $\tilde{a}_{3,2}=$0.118579912937683   \\
$\tilde{a}_{3,3}=$0.168770971875965    &$\tilde{a}_{4,1}=$ 0.172252519817939   & $\tilde{a}_{4,2}=$0.198870773989805  \\
$\tilde{a}_{4,3}=$ 0.011308123857846   & $\tilde{a}_{4,4}=$  0.121165193490219    & $\tilde{a}_{5,1}=$  0.107284797278235  \\
$\tilde{a}_{5,2}=$ 0.125275675715479   & $\tilde{a}_{5,3}=$ 0.003608216309487    & $\tilde{a}_{5,4}=$0.156483792310469  \\
$\tilde{a}_{5,5}=$ 0.069605617987546   & $\tilde{a}_{6,1}=$  0.111537098901469   & $\tilde{a}_{6,2}=$0.104025846350128   \\
$\tilde{a}_{6,3}=$0.228281462765514    & $\tilde{a}_{6,4}=$ 0.075203021889408   & $\tilde{a}_{6,5}=$0.214262575129975  \\ 
 $\tilde{a}_{7,1}=$ 0.115503702215039   &  $\tilde{a}_{7,2}=$0.116459412732323   &  $\tilde{a}_{7,3}=$0.152707824629209  \\
  $\tilde{a}_{7,4}=$ 0.080352146755342  &  $\tilde{a}_{7,5}=$ 0.242701834294582  &  $\tilde{a}_{7,6}=$ 0.261454998381366  \\
    \end{tabular}  
         
          \[    \mbox{and} \; \; \; \;    \vb = \vbt=   \left(   \begin{array}{l} 
          0.148802853943694\\
   0.140365832446254\\
   0.185913207665706\\
   0.143576841907452\\
   0.102077358038296\\
   0.109741290668591\\
   0.169522615330006   \end{array} \right).  \]

         \subsection{Methods for $K=100$}
          \subsubsection{An IMEX pair with $p^e=p^i=3$, $s=p_{lin}=5$}
             The nonzero coefficients in $\mA$ are given by: 
             
 \begin{tabular}{llll}
$a_{2,1}=$   0.607406844316321    & $a_{3,1}=$  0.330966515197897  & $a_{3,2}=$ 0.340310969038496  \\
 $a_{4,1}=$ 0.194835632796261    & $a_{4,2}=$ 0.050335014780643   & $a_{4,3}=$ 0.464427204928710      \\
$a_{5,1}=$ 0.135852828893193     & $a_{5,2}=$ 0.192467857403262   & $a_{5,3}=$ 0.024895163948772  \\
$a_{5,4}=$  0.337487088561988    & & \\
     \end{tabular}

         \begin{tabular}{llll}
 $\tilde{a}_{2,1}=$   0.607406844316321   & $\tilde{a}_{3,1}=$ 0.330966515197897 & $\tilde{a}_{3,2}=$   0.340310969038496 \\
 $\tilde{a}_{4,1}=$   0.193496010547777   &  $\tilde{a}_{4,2}=$0.200519538677067  &  $\tilde{a}_{4,3}=$ 0.088728444949044   \\
  $\tilde{a}_{4,4}=$0.226853858331728     &  $\tilde{a}_{5,1}=$ 0.129157547811257  &  $\tilde{a}_{5,2}=$ 0.131916477717161  \\
   $\tilde{a}_{5,3}=$ 0.231093457658500   &  $\tilde{a}_{5,4}=$0.037795421975484   &  $\tilde{a}_{5,5}=$0.160740033644814\\
        \end{tabular}  
        
  \[    \mbox{and} \; \; \; \;    \vb = \vbt=   \left(   \begin{array}{l} 
   0.247413560693329\\
   0.225966553626905\\
   0.158714688358981\\
   0.110694923985245\\
   0.257210273335540\\
   \end{array} \right).  \]

    \subsubsection{An IMEX pair with $s=10$, $p^e=4$, $p^i=3$, $p_{lin}=6$}
   
    The six stage IMEX pair with fourth order for nonlinear problems and sixth order on linear problems is given by the 
 nonzero coefficients in $\mA$: 
        
\begin{tabular}{llll}
$a_{2,1}=$ 0.296441642233457      &  $a_{3,1}=$0.296441642233457   & $a_{3,2}$= 0.296441642233457       \\
$a_{4,1}=$ 0.159129854065221      & $a_{4,2}=$0.159129854065221    & $a_{4,3}=$0.159129854065221    \\
$a_{5,1}=$0.019921058480607   & $a_{5,2}=$0.019307972560535   & $a_{5,3}=$0.019307972560535   \\
$a_{5,4}=$0.035968656715399   & $a_{6,1}=$0.083713252215858   & $a_{6,2}=$0.024819940994344   \\
$a_{6,3}=$0.014853700442704   & $a_{6,4}=$0.027670831336742   & $a_{6,5}=$0.228053739908418     \\
$a_{7,1}=$0.055091634074979   & $a_{7,2}=$0.016337375252376   & $a_{7,3}=$0.009777262948455 \\
$a_{7,4}=$0.018208562860852   & $a_{7,5}=$0.150068886917073  & $a_{7,6}=$ 0.195070983285382\\
$a_{8,1}=$0.049748358305482   &$a_{8,2}=$0.025038891947022   &$a_{8,3}=$0.020856206502734   \\
$a_{8,4}=$0.017683704471081   & $a_{8,5}=$0.095682957831884   &$a_{8,6}=$0.124376005255730   \\
$a_{8,7}=$0.189009285909559   &$a_{9,1}=$0.073381984734260  &$a_{9,2}=$0.048836937556200\\
$a_{9,3}=$0.030836020427511 &$a_{9,4}=$0.039354142711587   & $a_{9,5}=$0.082928693754279  \\
$a_{9,6}=$0.106285066647655   & $a_{9,7}=$0.161517203488072  & $a_{9,8}=$0.253323136059411  \\
$a_{10,1}=$0.065586205155812 &  $a_{10,2}=$0.044449329946329& $a_{10,3}=$0.028936784917608   \\
$a_{10,4}=$0.033888353855577  & $a_{10,5}=$0.071414307430042  & $a_{10,6}=$0.091527804079801 \\
$a_{10,7}=$0.139082279606173  & $a_{10,8}=$0.218136263377849   &$a_{10,9}=$0.255265559839090 \\
  \end{tabular}

   \noindent The nonzero coefficients in $\mAt$ are given by: 
   
   \begin{tabular}{llll}
$\tilde{a}_{2,1}=$-0.047916334864055   & $\tilde{a}_{2,2}=$ 0.344357977097512   &$\tilde{a}_{3,1}=$-0.082116017384178 \\
$\tilde{a}_{3,2}=$ 0.330641324753580   & $\tilde{a}_{3,3}=$ 0.344357977097512 & $\tilde{a}_{4,1}=$0.948951469602658 \\
$\tilde{a}_{4,2}=$-0.238288155973038   & $\tilde{a}_{4,3}=$-0.577631728531468 & $\tilde{a}_{4,4}=$0.344357977097512 \\
$\tilde{a}_{5,1}=$0.441387001271949    & $\tilde{a}_{5,2}=$-0.163364374457318  & $\tilde{a}_{5,3}=$-0.504386376049355 \\
$\tilde{a}_{5,4}=$-0.023488567545713   & $\tilde{a}_{5,5}=$0.344357977097512   & $\tilde{a}_{6,1}=$ 0.099545320449493   \\
$\tilde{a}_{6,2}=$1.759501344125055    & $\tilde{a}_{6,3}=$-0.824006411866389  & $\tilde{a}_{6,4}=$1.409963103927893 \\
$\tilde{a}_{6,5}=$-2.410249868835498   &$\tilde{a}_{6,6}=$0.344357977097512    &    $\tilde{a}_{7,1}=$0.151486935292728 \\
$\tilde{a}_{7,2}=$0.407347161353228    & $\tilde{a}_{7,3}=$-0.244112999903862  & $\tilde{a}_{7,4}=$0.557935535998744 \\
$\tilde{a}_{7,5}=$-0.774009241132186   & $\tilde{a}_{7,6}=$0.001549336632953   & $\tilde{a}_{7,7}=$ 0.344357977097512 \\
$\tilde{a}_{8,1}=$1.113607828755682   & $\tilde{a}_{8,2}=$ -1.940412832219541  & $\tilde{a}_{8,3}=$1.777801724820027 \\
$\tilde{a}_{8,4}=$ -0.950842851699990 &$\tilde{a}_{8,5}=$1.780109034747002   & $\tilde{a}_{8,6}=$0.158989846049901 \\
$\tilde{a}_{8,7}=$-1.761215317327101   & $\tilde{a}_{8,8}=$0.344357977097512  &  $\tilde{a}_{9,1}=$-1.089122782612700 \\
$\tilde{a}_{9,2}=$0.803069401235045    & $\tilde{a}_{9,3}=$-0.090413146702761  & $\tilde{a}_{9,4}=$-0.697564076196375 \\
$\tilde{a}_{9,5}=$0.822000558154658    & $\tilde{a}_{9,6}=$-0.611804180680986   & $\tilde{a}_{9,7}=$1.687933136273893 \\
$\tilde{a}_{9,8}=$-0.371993701189312   & $\tilde{a}_{9,9}=$0.344357977097512    &  $\tilde{a}_{10,1}=$0.444596175510804 \\
$\tilde{a}_{10,2}=$-0.716394629436675 & $\tilde{a}_{10,3}=$0.093268024454960  & $\tilde{a}_{10,4}=$0.501655569949843 \\
$\tilde{a}_{10,5}=$0.351539276646713   & $\tilde{a}_{10,6}=$0.952014224936281  & $\tilde{a}_{10,7}=$-0.598231110418698 \\
$\tilde{a}_{10,8}=$0.115776872278219   & $\tilde{a}_{10,9}=$-0.540295492810681  & $\tilde{a}_{10,10}=$0.344357977097512\\
  \end{tabular}  
   
   \smallskip

\[ \mbox{and} \; \; \; \; \vb = \vbt=  \left(    \begin{array}{l}  
0.084877037374285 \\ 
0.060614911672550 \\
0.050348209787848\\
0.029884478272179\\
0.088565247859017 \\
0.114451890708771\\
0.173923410482789 \\
0.112585396241353 \\
0.131748723280703 \\
0.153000694320504 \\
    \end{array}  \right) \]

\bibliography{LNL}
\end{document}